\documentclass[oneside,a4paper]{amsart} % Temporarily "oneside" to have todonotes always on the right
\usepackage{amssymb,verbatim,mathabx}

\usepackage[shortlabels]{enumitem}
    \setenumerate[0]{label={\rm (\roman*)}, leftmargin=*}

\usepackage[T1]{fontenc}
\usepackage[english]{babel}

\usepackage[textsize=footnotesize,textwidth=20ex,colorinlistoftodos]{todonotes}

\usepackage{hyperref}
\hypersetup{
    colorlinks=true,
    linkcolor=blue,
    filecolor=magenta, 
    citecolor=red,     
    urlcolor=cyan,
}

\usepackage[capitalize]{cleveref}
    \crefname{enumi}{}{}
    \Crefname{enumi}{Item}{Items}
    \crefname{equation}{}{}
    \Crefname{equation}{Equation}{Equations}

\newtheorem{proposition}{Proposition}[section]
\newtheorem{lemma}[proposition]{Lemma}
\newtheorem{corollary}[proposition]{Corollary}
\newtheorem{theorem}[proposition]{Theorem}

\theoremstyle{definition}

\newtheorem{example}[proposition]{Example}

\newtheorem{remark}[proposition]{Remark}

\newtheorem{algorithm}[proposition]{Algorithm}

\newcommand{\Rad}{{\rm Rad}}

\newcommand{\one}{\mathbf{1}}

\begin{document}
\title{Quasi-Clifford Algebras, their Lie and Jordan Algebras, and Quadratic Forms over $\mathbb{F}_2$}
\author{Hans Cuypers}

\begin{abstract}
Let $R$ be a commutative and associative ring  containing distinct elements $1$ and $-1$. 
Let $\Gamma=(\mathcal{V},\mathcal{E},\lambda)$ be a labeled graph, with vertex set $\mathcal{V}$, edge set $\mathcal{E}$ and labeling of the vertices $\lambda:\mathcal{V}\rightarrow R^*$ which maps $v\in \mathcal{V}$ to invertible elements $\lambda(v)\in R^*$.
 Then we consider  the associative $R$-algebra $\mathfrak{C}(\Gamma)$ with identity element $\one$ generated by the elements of  $\mathcal{V}$ such that for all
 $v,w\in \mathcal{V}$ we have
 
$$\begin{array}{lll} 
  v^2 &=\lambda(v)\one,&\\
  vw+wv&=0&\textrm{if } \{v,w\}\in \mathcal{E},\\
  vw-wv&=0&\textrm{if } \{v,w\}\not\in \mathcal{E}.\\
\end{array}$$

If $\Gamma$ is the complete graph, $\mathfrak{C}(\Gamma)$ is a Clifford algebra, otherwise it is a so-called quasi-Clifford algebra.

We describe this algebra as a twisted group algebra with the help of a    vector space $V$ over the field $\mathbb{F}_2$
equipped with a bilinear form $g$.
Using this description, we  determine the isomorphism type of $\mathfrak{C}(\Gamma)$ for several interesting graphs $\Gamma$.

As the algebra $\mathfrak{C}(\Gamma)$ is associative, we can also consider the corresponding Lie algebra with Lie bracket $[\cdot,\cdot]$ and Jordan algebra with multiplication $\circ$, as well as  some of their subalgebras.
We find that the elements  $v,w\in \mathcal{V}$ satisfy the following relations
$$\begin{array}{lll} 
  {[v,w]}&=0&\textrm{if } \{v,w\}\not\in \mathcal{E},\\
  {[v,[v,w]]}&=\lambda(v) w&\textrm{if } \{v,w\}\in \mathcal{E}.\\
\end{array}$$
and 
$$\begin{array}{lll}
     v\circ v&=\lambda(v)\one\\
  {v\circ w}&=0&\textrm{if } \{v,w\}\in \mathcal{E},\\
  {v\circ(v\circ w)}&=\lambda(v) w&\textrm{if } \{v,w\}\not\in \mathcal{E}.\\
\end{array}$$

We provide characterizations of both the Lie and Jordan algebras generated by the elements in $\mathcal{V}$, as algebras defined by these relations.

In case $R$ is a field of characteristic $0$, we can identify these Lie algebras with quotients of the compact 
subalgebras of Kac-Moody Lie algebras and prove that they admit a so-called generalized spin representation. 
\end{abstract}

\maketitle

\section{Introduction}

In this paper $R$ is an  associative and commutative ring containing distinct elements $1$ and $-1$ in which $2$ is not a zero divisor.
In case $R$ is a field (of characteristic different from $2$), we also denote it by $\mathbb{F}$.
We will study certain algebras over $R$, in particular Lie and Jordan algebras, which are defined by generators and relations.
These generators and relations are defined by labeled graphs.

So, let $\Gamma=(\mathcal{V},\mathcal{E},\lambda)$ be a labeled graph with vertex set $\mathcal{V}$, edge set $\mathcal{E}$ and labeling $\lambda:\mathcal{V}\rightarrow R^*$,  the set of invertible elements from $R$, where a vertex $v\in\mathcal{V}$ is  labeled by the  invertible element $\lambda(v)$.   
Then, before considering Lie and Jordan algebras,  we first consider the associative algebra $\mathfrak{C}(\Gamma)$ over $R$ with identity element $\one$ generated 
by the elements of  $\mathcal{V}$ such that for all
 $v,w\in \mathcal{V}$ we have
 $$\begin{array}{lll} 
  v^2 &=\lambda(v)\one,&\\
  vw+wv&=0&\textrm{if } v\sim w,\\
  vw-wv&=0&\textrm{if } v\not \sim w.\\
\end{array}$$
Here $v\sim w$ denotes that $\{v,w\}$ is an edge in $\mathcal{E}$.

If $\Gamma$ contains no edges and  $\lambda(v)=1$ for all $v\in \mathcal{V}$, then the algebra $\mathfrak{C}(\Gamma)$ is a Grassmann algebra.
On the other hand, if $\Gamma$ is the complete graph on $n$ vertices, $R=\mathbb{R}$ and $\lambda(v)=\pm 1$ for all $v\in \mathcal{V}$, then the algebra $\mathfrak{C}(\Gamma)$ is a real Clifford algebra $\mathrm{Cl}(p,q)$,
where $n=p+q$ and $p$ vertices are labeled with $+1$, while $q$ vertices have the label $-1$.

For an arbitrary finite graph $\Gamma$ and field $R$,  we obtain a so-called quasi-Clifford algebra as studied by Gastineau-Hills in \cite{quasi} in connection with orthogonal designs (see also \cite{leopardi1,Leopardi2,Seberry1} and the recent book \cite{Seberry}).

In this paper we  describe for arbitrary such labeled graphs $\Gamma$ the algebra $\mathfrak{C}(\Gamma)$ as a twisted group algebra
$\mathfrak{C}(V,g,\lambda)$.
Here $V=V_\Gamma$ is the $\mathbb{F}_2$-space  with basis $\mathcal{V}$ and addition $\pluscirc$ equipped with a bilinear form $g$, and
map $\lambda:\mathcal{V}\rightarrow R^*$  which assigns to each element $v$ in $\mathcal{V}$ the invertible element
$\lambda(v)$ from $R$. The algebra $\mathfrak{C}(V,g,\lambda)$ is an $R$ algebra with basis $V$.
Multiplication is  defined as the bilinear expansion of the product  defined on $v,w\in V$ by
$$vw=(-1)^{g(v,w)}\prod_{u\in \mathcal{V},\\ u^*(v)u^*(w)=1}\ (-1)^{g(u,u)} \lambda(u)\ v\pluscirc w,$$
where for $u\in \mathcal{V}$ the element $u^*\in V^*$ satisfies $u^*(v)=\delta_{u,v}$ for all  $v\in\mathcal{V}$.

Notice that the vector space $V$  is equipped with the alternating form $f=f_\Gamma$ defined by $f(v,w)=1$ for $v,w\in \mathcal{V}$ if and only if $v$ and $w$ are adjacent.
In particular, we prove:

\begin{theorem}\label{mainclifford}
Let $\Gamma=(\mathcal{V},\mathcal{E},\lambda)$ be a labeled graph.
Then, for every bilinear form  $g$ on $V$ with $g(v,w)+g(w,v)=f_\Gamma(v,w)$  for $v, w\in \mathcal{V}$, we find  $\mathfrak{C}(\Gamma)$ to be isomorphic to the $R$-algebra $\mathfrak{C}(V,g,\lambda)$.
\end{theorem}

The proof of this result is given in the Sections \ref{sect:twisted} and \ref{sect:generators}.

Notice that the above theorem implies that the algebras $\mathfrak{C}(V,g,\lambda)$ and $\mathfrak{C}(V,g',\lambda)$
where both $g$ and $g'$ are bilinear forms
with  $g(v,w)+g(w,v)=g'(v,w)+g'(w,v)=f_\Gamma(v,w)$ for $v, w\in \mathcal{V}$ are isomorphic.

This theorem extends some of the results of \cite{Ovsienko,Twisted_clifford} and \cite{Elduque_clifford}, which focus on
finite dimensional Clifford algebras (and hence on finite complete graphs $\Gamma$ with labels $\pm 1$).
Although the proof of \cref{mainclifford} follows more or less the line of \cite{Ovsienko,Twisted_clifford,Elduque_clifford}, 
we provide it in the Sections
\ref{sect:twisted} and \ref{sect:generators}, as the notions introduced in these sections are of great use in the remaining sections.

Given the description as twisted group algebra we determine the structure of the algebra $\mathfrak{C}(\Gamma)$ in \cref{sect:structure}, and  determine its isomorphism type  for several interesting graphs $\Gamma$ with labels $\pm 1$. (Notice that we can assume that labels are $1$ in case $R$ is an algebraically closed field, or $\pm 1$ in case $R=\mathbb{R}$, as we can replace generators by scalar multiples.) This is done in Section \ref{sect:examples}. 
As an example we apply our results to complete graphs and obtain quickly the classification of Clifford algebras. 
But we also consider graphs of type $A_n, D_n$ and $E_n$, as they are related to Lie algebras studied in later sections.

Indeed, as the algebra $\mathfrak{C}:=\mathfrak{C}(\Gamma)$ is associative, we can also consider the corresponding Lie algebra $\mathfrak{C}^-$ with Lie product $[\cdot,\cdot]$ defined by 
$$[x,y]=\frac{1}{2}(xy-yx),$$
as well as the Jordan algebra $\mathfrak{C}^+$ with Jordan product $\circ$ defined by 
$$x\circ y=\frac{1}{2}(xy+yx)$$
and some of their subalgebras.  (Here we assume $R$ to contain the element $\frac{1}{2}$.)
In particular, we investigate the Lie and Jordan  subalgebras of $\mathfrak{C}^-$ and $\mathfrak{C}^+$, respectively, 
generated by the elements of $\mathcal{V}$. 

The elements  $v\neq w\in \mathcal{V}$ satisfy the following relations:
$$\begin{array}{lll} 
   v^2&=\lambda(v)\one,&\\
  {[v,w]}&=0&\textrm{if } v\not\sim w,\\
  {[v,[v,w]]}&= \lambda(v) w&\textrm{if } v \sim w.\\
\end{array}$$

For the  Jordan product $\circ$ we have
$$\begin{array}{lll} 
  v\circ v&=\lambda(v)\one,&\\
  v\circ w&=0&\textrm{if } v\sim w,\\
  v\circ(v\circ w)&=\lambda(v) w&\textrm{if } v\not \sim w.\\
\end{array}$$

We provide characterisations of the subalgebras of $\mathfrak{C}^-$ and   $\mathfrak{C}^+$ generated by the elements of $\mathcal{V}$, as the universal Lie or Jordan algebras satisfying the above relations and the additional requirement that they are embedded in the Lie or Jordan algebra of an associative algebra.

As indicated above, the algebra $\mathfrak{C}=\mathfrak{C}(\Gamma)$ can be described as a twisted group algebra with the help of an $\mathbb{F}_2$-vector space $V=V_\Gamma$ equipped with a bilinear form $g$. 
Let $Q$ be a  quadratic form   on $V$ satisfying
$$Q(v+w)+Q(v)+Q(w)=g(v,w)+g(w,v)$$
for all $v,w\in V$.
We say $Q$ is \emph{compatible} with $\Gamma$.
Then we can define the involution $$\tau_Q:\mathfrak{C}\rightarrow \mathfrak{C}$$
which maps all elements of $V$ with $Q(v)=1$ to $-v$ and fixes all others.
This involution $\tau_Q$ is an anti-automorphism of $\mathfrak{C}(\Gamma)$. 
See \cref{sect:involutions}.

In the Sections \ref{sect:lie}, \ref{sect:liegraphs} and \ref{sect:kac} we then use the geometry of the space $V$
induced by the alternating form $f$, defined by $f(u,v)=g(u,v)+g(v,u)$ for all $u,v\in V$, 
to classify the Lie subalgebras of $\mathfrak{C}^-$ generated the elements of $\mathcal{V}$ and prove the following result for Lie algebras.

\begin{theorem}\label{mainlie}
Assume $R$ contains $\frac{1}{2}$ and
let $\Gamma=(\mathcal{V},\mathcal{E},\lambda)$ be a labeled connected graph with labels in $R^*$ and $\mathfrak{C}=\mathfrak{C}(\Gamma)$.

Suppose $\mathfrak{A}$ is an associative algebra containing  elements $\overline{v}$, with $v\in\mathcal{V}$, satisfying the following relations 
in $\mathfrak{A}^-$:
$$\begin{array}{lll}
   \overline{v}^2&= \lambda(v)\one_{\mathfrak{A}} &  \\
   {[\overline{v},\overline{w}]}&  =0 & \textrm{if } v\not\sim w,\\
   {[\overline{v},[\overline{v},\overline{w}]]}&= \lambda(v) \overline{w} & \textrm{if } v \sim w.\\
\end{array}$$

Then, with $Q$ being  the unique quadratic form  compatible with $\Gamma$ and with $Q(v)=1$ for all $v\in \mathcal{V}$, we find one of the following:
\begin{enumerate}
\item The subalgebra of  $\mathfrak{A}^-$ generated by the elements $\overline{v}$, with $v\in\mathcal{V}$, is isomorphic to a quotient of the Lie subalgebra $\{x\in \mathfrak{C}^-\mid \tau_Q(x)=-x\}$ of $\mathfrak{C}^-$ of skew elements with respect to the involution $\tau_{Q}$.
\item There are two disjoint  sets $\Omega,\Omega'$, such that the  subalgebra of $\mathfrak{A}^-$ generated by the elements $\overline{v}$, with $v\in\mathcal{V}$, is isomorphic to a quotient of the Lie algebra $\mathfrak{g}(\Omega,\Omega')$, which is contained in the Lie subalgebra of $\mathfrak{C}^-$  of skew elements with respect to the involution $\tau_{Q}$.
\end{enumerate}
\end{theorem}

(For definition and notation see Section \ref{sect:liegraphs}.)

We notice that the second case of the conclusion of the theorem occurs when the graph $\Gamma$ is  a line graph of a multi-graph, i.e. a graph whose vertices are the edges of a multi-graph $\Delta$, and two such edges adjacent if and only if they meet in a single vertex of $\Delta$.

These Lie algebras also appear in \cite{Tanya1,gintz}, where  ordinary finite graphs and $R$ the field of complex numbers are considered.
In case $R$ is a field of characteristic $0$ and $\lambda$ takes only the value $-1$, we can by \cite{Berman} identify these Lie algebras with quotients of  compact 
subalgebras of Kac-Moody Lie algebras and prove that they admit a so-called generalized spin representation.
In particular, using the computations of Section \ref{sect:examples} and \ref{sect:lie}, we are able to identify various quotients of these 
compact Lie subalgebras of Kac-Moody algebras and construct spin representations of such algebras extending the results of \cite{ Damour,Buyl,spin}. This is the topic of Section \ref{sect:kac}.

For Jordan algebras we obtain a similar characterization. We prove:

\begin{theorem}\label{mainjordan}
Assume $\frac{1}{2}\in R$ and let $\Gamma=(\mathcal{V},\mathcal{E},\lambda)$ be a connected labeled graph with labels in $R^*$, with connected complement.
Let $\mathfrak{C}$ be the algebra $\mathfrak{C}(\Gamma)$.

Suppose $\mathfrak{A}$ is an associative algebra containing  elements $\overline{v}$, with $v\in\mathcal{V}$, satisfying the following relations 
in $\mathfrak{A}^+$:
$$\begin{array}{lll}
   \overline{v}\circ \overline{v}&= \lambda(v)\one_{\mathfrak{A}} &  \\
   {\overline{v}\circ\overline{w}}&  =0 & \textrm{if } v\sim w,\\
   {\overline{v}\circ(\overline{v}\circ\overline{w})}&= \lambda(v) \overline{w} & \textrm{if } v \not\sim w.\\
\end{array}$$
Here $\circ$ denotes the Jordan product.

Then the subalgebra of  $\mathfrak{A}^+$ generated by the elements $\overline{v}$, with $v\in\mathcal{V}$, is isomorphic to 
a quotient of the subalgebra $\{x\in \mathfrak{C}^+\mid \tau_Q(x)=x\}$ of $\mathfrak{C}^+$ of symmetric elements with respect to the involution $\tau_{Q}$, where  $Q$ is the unique quadratic form compatible with $\Gamma$ and with $Q(v)=0$ for all $v\in \mathcal{V}$.
\end{theorem}

A proof of this result can be found in \cref{sect:jordan}.

\section{A class of algebras obtained from bilinear forms over $\mathbb{F}_2$}

\label{sect:twisted}

In this section we provide a description of a class of  twisted group algebras. In the case that they are finite dimensional, these  algebras  turn out to be  quasi-Clifford algebras as introduced by Gastineau-Hills \cite{quasi}. Our description as twisted group algebra is closely related to the description of Clifford algebras as twisted group algebras, see \cite{Twisted_clifford}, and relates our algebras
to quadratic spaces over the field with two elements as in \cite{Ovsienko,Elduque_clifford}. (See also the work of Shaw \cite{shaw1,shaw2,shaw3,shaw4}.)

We start with an $\mathbb{F}_2$ vector space $V$ (with addition $\pluscirc$) equipped with a bilinear form $g:V\times V\rightarrow \mathbb{F}_2$.
Let $\mathcal{V}$ be a basis for $V$ and $\mathcal{V}^*$ a dual basis, where $v^*$ denotes the dual of $v\in \mathcal{V}$,
i.e. $v^*(v)=1$ and $v^*(w)=0$ for all $w\in \mathcal{V}\setminus \{v\}$.
As in the introduction, let  $R$  denote a commutative and associative ring, and $R^*$  its set of invertible elements including the distinct elements $1$ and $-1$.
Finally, let $\lambda:\mathcal{V}\rightarrow R^*$ be a map which we extend to
a map $\Lambda: V\times V\rightarrow R^*$ by 

$$\Lambda(v,w):=\prod_{b\in \mathcal{V},\\ b^*(v)b^*(w)=1}\ (-1)^{g(b,b)}\lambda(b)$$
for all $v,w\in V$. 

Notice that this is well defined, also for infinite dimensional spaces $V$, since almost all values of $b^*(v)b^*(w)$, where $b\in \mathcal{V}$, are $0$. Moreover, for all $u\in V$ we have $\Lambda(u,0)=\Lambda(0,u)=1$.

The twisted group algebra $\mathfrak{C}{(V,g,\lambda)}$ is then the $R$-algebra with
basis $\{v\mid v\in V\}$, unit element $0\in V$ which we denote by $\one$,  and
bilinear multiplication  defined by

$$
\begin{array}{ll}
  vw&=(-1)^{g(v,w)}\cdot\Lambda(v,w)\ {v\pluscirc w}\\
\end{array}$$
for all $v,w\in V$.

Following Gastineau-Hills \cite{quasi}, we will call these algebras \emph{quasi-Clifford algebras}.

Notice that elements $v$ and $w$, where $v\neq w\in V$, satisfy the relations
$$\begin{array}{lll}
  vw-wv=0& \textrm{if } f(v,w)=0\\
   vw+wv=0& \textrm{if } f(v,w)=1\\
\end{array}$$
for all $v,w\in V$, where $f$ is the bilinear (and also alternating) form given by $f(v,w)=g(v,w)+g(w,v)$.
%The algebra $\mathfrak{C}{(V,Q,\lambda)}$ is a twisted group algebra as considered in \cite{Elduque_clifford,Twisted_clifford}.
Moreover, for $v\in \mathcal{V}$ we have 
$$v^2=\lambda(v)\one.$$

\begin{proposition}
The algebra $\mathfrak{C}{(V,g,\lambda)}$ is  associative.
\end{proposition}

\begin{proof}
  Let $u,v,w\in V$, then
  $$\begin{array}{ll}
    u(vw)&=u\left((-1)^{g(v,w)}\Lambda(v,w)\ {v\pluscirc w}\right)\\
    &=(-1)^{g(v,w)+g(u,v\pluscirc w)}\Lambda(u,v\pluscirc w)\Lambda(v,w)\ {u\pluscirc v\pluscirc w}\\
    &=(-1)^{g(v,w)+g(u,v)+g(u,w)}\Lambda(u,v\pluscirc w)\Lambda(v,w)\ {u\pluscirc v\pluscirc w},
    \end{array}$$
  while
  $$\begin{array}{ll}
    (uv)w&=\left((-1)^{g(u,v)}\Lambda(u,v)\ {u\pluscirc v}\right)w\\
    &=(-1)^{g(u,v)+g(u\pluscirc v,w)}\Lambda(u,v)\Lambda(u\pluscirc v,w)\ {u\pluscirc v\pluscirc w}\\
&  =(-1)^{g(u,v)+g(u,w)+g(v,w)}\Lambda(u,v)\Lambda(u\pluscirc v,w)\ {u\pluscirc v\pluscirc w}.\end{array}$$  

  So, we find the algebra to be associative,
  if and only if $\Lambda$ satisfies
  $$\Lambda(u,v)\cdot\Lambda (u\pluscirc v,w)=\Lambda(v,w)\cdot \Lambda(u,v\pluscirc w)$$
  for all $u,v,w\in V$.
  
  This latter equality is obvious if $u=w$. If $u=v$, then this equation reads  as  $\Lambda(u,u)=\Lambda(u,w)\Lambda(u,u\pluscirc w)$
  and follows from the observation that for $b\in \mathcal{V}$ we have   $b^*(u)=1$ if and only if $b^*(w)=1$ and $b^*(u\pluscirc w)=0$
  or  $b^*(w)=0$ and $b^*(u\pluscirc w)=1$.
  Similarly we find that for $v=w$ the identity holds.
  For distinct  $u,v,w$ the identity follows from the observation that for all
  $b^*\in V^*$ we have
  $$\begin{array}{lll}
    b^*(u)b^*(v)=1\Leftrightarrow b^*(u\pluscirc v)=0&\mathrm{and}&b^*(v)b^*(w)=1\Leftrightarrow b^*(v\pluscirc w)=0.\\
    \end{array}$$
    Indeed, it follows that both the left hand side as well as the right hand side of the equation is the product of those terms $(-1)^{g(b,b)}\lambda(b)$, where $b\in \mathcal{V}$ for which at least two of the three values $b^*(u)$, $b^*(v)$ and $b^*(w)$ are $1$.
  \end{proof}

\begin{remark}\label{changeofbasis}
The definition of the algebra $\mathfrak{C}(V,g,\lambda)$ depends on $V$, $g$ and $\lambda:\mathcal{V}\rightarrow R^*$, and hence implicitly on the basis $\mathcal{V}$ of $V$.

In the next section we will see that replacing $\mathcal{V}$ by a new basis $\mathcal{V}'$ and $\lambda$ by $\lambda':\mathcal{V}'\rightarrow R^*$ where
$\lambda(v')=(-1)^{g(v',v')}\Lambda(v',v')$ for $v'\in \mathcal{V}'$ does not change the isomorphism type of the algebra.

 \end{remark}
 
\section{From relations to algebra}
\label{sect:generators}
\label{sect:relations}

Let $R$ be a commutative and associative ring with distinct elements $1,-1$.
Suppose $V$ is an $\mathbb{F}_2$-space equipped with a bilinear form $g$ and for some basis $\mathcal{V}$ of $V$ a map $\lambda:\mathcal{V}\rightarrow R^*$ which we extend to a map 
 $\Lambda:V\times V\rightarrow R^*$ defined by $\Lambda(v,w):=\prod_{b\in \mathcal{V},b^*(v)b^*(w)=1}\ (-1)^{g(b,b)} \lambda(b)$ for all $v,w\in V$.
 Then we can consider the algebra $\mathfrak{C}(V,g,\lambda)$ over $R$ as defined in \cref{sect:twisted}.
 
The elements $v\neq w\in \mathcal{V}$ satisfy the following relations in $\mathfrak{C}(V,g,\lambda)$:
 $$
 \begin{array}{ll}
 v^2&=\lambda(v)\one\\
 vw&=(-1)^{f(v,w)}wv,\\
 \end{array}$$
 for all $v,w\in V$, where $f$ is the alternating form defined by $f(v,w)=g(v,w)+g(w,v)$.
 
 We can capture this information in a labeled graph. This graph has vertex set  $\mathcal{V}$, two vertices
$v\neq w$ being adjacent if and only if $vw=-wv$. A vertex $v\in\mathcal{V}$ is labeled by $\lambda(v)$ if and only if $v^2=\lambda(v) \one$.

In this section we reverse the above process by showing that each such labeled graph determines  an associative algebra, given by  generators and relations, isomorphic to an algebra  $\mathfrak{C}(V,g,\lambda)$.

So, let $\Gamma=(\mathcal{V},\mathcal{E},\lambda)$ be a labeled graph 
with vertex set $\mathcal{V}$ and edge set $\mathcal{E}$,
and the vertices $v\in \mathcal{V}$
labeled by nonzero invertible elements $\lambda(v)$
from  $R$. 

Denote by $\mathfrak{C}(\Gamma)$ the 
associative algebra generated by the elements from $\mathcal{V}$ subject to the relations

$$\begin{array}{lll} 
  v^2 &=\lambda(v)\one,&\\
  vw+wv&=0&\textrm{if } v\sim w,\\
  vw-wv&=0&\textrm{if } v\not \sim w.\\
\end{array}$$

Put a linear  order $<$ on $\mathcal{V}$ and consider ${V}_\Gamma$, the vector space of finite ordered
subsets of $\mathcal{V}$, where for two finite subsets $v,w$ of $\mathcal{V}$ the sum
$v\pluscirc w$ is defined to be the ordered symmetric difference of $v$ and $w$.
Then the elements of $\mathfrak{C}(\Gamma)$ can be identified with 
$R$-linear combinations of elements from $V_\Gamma$.

The product of $\mathfrak{C}(\Gamma)$ can then be defined in the following way.

For $v$ and $w$ being two finite subsets of $\mathcal{V}$, let $g_\Gamma(v,w)$ denote the number of pairs $(x,y)\in v\times w$, where  
$\{x,y\}$ is an edge with $x<y$, modulo $2$.
Then  $g_\Gamma(u,v\pluscirc w)=g_\Gamma(u,v)+g_\Gamma(u,w)$, for any finite subsets $u,v,w$ of $\mathcal{V}$,
as the  edges $\{x,z\}$ with $x<z$,  $x\in u$ and  $z\in v\cap w $, are counted twice at the right hand site of the equation.
Similarly $g_\Gamma(v\pluscirc w,u)=g_\Gamma(v,u)+g_\Gamma(w,u)$, for any finite subsets $u,v,w$ of $\mathcal{V}$.
So, $g_\Gamma$ is a bilinear map from $V$ to $\mathbb{F}_2$.
For $v\in \mathcal{V}$ we have $g_\Gamma(v,v)=0$.

The product of  $v$ and $w$ in $\mathfrak{C}(\Gamma)$ is then 
$$vw=(-1)^{g_\Gamma(v,w)}\left(\prod_{x\in v\cap w}\lambda(x)\right)\, v\pluscirc w.$$

Clearly this definition of the product is forced upon us by the relations and associativity of the product.
The sign $(-1)^{g_\Gamma(v,w)}$ relates to the ordering of the elements and the relations
for two adjacent vertices of the graph, while the product  $\prod_{x\in v\cap w}\lambda(x)$ relates to the labeling of the vertices.

By construction, the algebra $\mathfrak{C}(\Gamma)$ is the universal associative algebra satisfying the relations prescribed by the graph $\Gamma$.
So, we have proven the following.

\begin{theorem}\label{graphalgebra}
The algebra $\mathfrak{C}(\Gamma)$ is isomorphic to $\mathfrak{C}(V_\Gamma,g_\Gamma,\lambda)$.
\end{theorem}

Moreover, we find:

\begin{theorem}\label{isomorphic}
Let $V$ be an $\mathbb{F}_2$-space with basis $\mathcal{V}$, equipped with bilinear forms $g$ and $g'$  and map $\lambda:\mathcal{V}\rightarrow R^*$. 

Then the algebra $\mathfrak{C}(V,g,\lambda)$ is isomorphic
to $\mathfrak{C}(\Gamma)$ where  
$\Gamma$ is the graph with vertex set $\mathcal{V}$, in which two vertices $v,w$ are adjacent if and only if
$f(v,w)=g(v,w)+g(w,v)=1$,  and a  vertex $v$  is labeled
by $\lambda(v)$ if and only if $v^2=\lambda(v)\one$.

In particular, $\mathfrak{C}(V,g,\lambda)$ and $\mathfrak{C}(V,g',\lambda)$ are isomorphic if
$g(u,v)+g(v,u)=g'(u,v)+g'(v,u)$ for all $u,v\in \mathcal{V}$.
\end{theorem}

The above theorems imply \cref{mainclifford}.

We also have the following consequence of the above theorems.
As before, let $V$ be an $\mathbb{F}_2$-space with basis $\mathcal{V}$ and equipped with bilinear form $g$ and map $\lambda:\mathcal{V}\rightarrow R^*$
and consider the algebra $\mathfrak{C}=\mathfrak{C}(V,g,\lambda)$.

Suppose $\mathcal{V'}$ is a second basis for $V$, and $\lambda':\mathcal{V}'\rightarrow R^*$ is given by
$\lambda'(v')=(-1)^{g(v',v')}\Lambda(v',v')$ for all $v'\in \mathcal{V}'$, where $\Lambda$ is defined as in the previous section.
Let $\Gamma'$ be the labeled graph with vertex set $\mathcal{V}'$, two vertices $v'$ and $w'$ adjacent if and only if $v'w'+w'v'=0$
and labeled by $\lambda'$.

Then the map that takes each vertex from $\mathcal{V'}$ to the corresponding element of $\mathfrak{C}$
induces a surjective homomorphism from $\mathfrak{C}(\Gamma')$ to $\mathfrak{C}$.
This homomorphism
is actually an isomorphism. Indeed, an element $k$ of the kernel of this map is contained  in a subalgebra of $\mathfrak{C}(\Gamma')$
generated by a finite subset $\mathcal{V}'_0$ of $\mathcal{V'}$. This subalgebra has dimension $2^{|\mathcal{V}'_0|}$ and maps 
surjectively onto the subalgebra of $\mathfrak{C}$ generated by $\mathcal{V}'_0$, which has the same dimension.
But then this homomorphism restricted to the subalgebras is an isomorphism, proving  
$k$ to be $0$.

So, we find the following:

\begin{corollary}\label{differentgraphs}
$\mathfrak{C}(\Gamma')$ is isomorphic to $\mathfrak{C}(V,g,\lambda)$.
\end{corollary}

This result will be used without reference in the sequel.

\section{Special quasi-Clifford algebras}
\label{sect:special}

The  case that $\Gamma=(\mathcal{V},\mathcal{E},\lambda)$ is a labeled graph with all labels $\pm 1$ deserves special attention.
Indeed, in the case that $R$ is an algebraically closed field we can assume that all labels are $1$. If $R$ is $\mathbb{Z}$, $\mathbb{R}$
or a finite field of odd order in which $- 1$ is not a square, all labels can assumed to be $\pm 1$.
So suppose all labels of $\Gamma$ are $\pm 1$.

Let $V_\Gamma$ the $\mathbb{F}_2$-space with basis $\mathcal{V}$ equipped with bilinear form $g_\Gamma$ as in the previous section.
By $f_\Gamma$ we denote the symplectic form defined by $$f_\Gamma(v,w)=g_\Gamma(v,w)+g_\Gamma(w,v)$$
for all $v,w\in V$.

Then let  $g$ be a (second) bilinear form on $V=V_\Gamma$
with $$g(u,v)+g(v,u)=f_\Gamma(u,v)$$ 
for all $u\neq v\in \mathcal{V}$ and additionally 
$$g(v,v)=\begin{cases} 1& \mathrm{if}\ \lambda(v)=-1,\\
0& \mathrm{if}\
\lambda(v)=1.\\
\end{cases}$$

Such a form does exist.
For example, suppose $<$ is a linear ordering on $\mathcal{V}$, then for $u$ and $w$ in $V$ we can define  $g(u,w)$ to be equal to the number of ordered pairs $(x,y)\in u\times w$, where $x<y$ and $\{x,y\}$ is an edge, or $x=y$ has label $-1$, modulo $2$.

Multiplication in $\mathfrak{C}(\Gamma)$ is then
define by
$$vw=(-1)^{g(v,w)}v\pluscirc w,$$
for all $v,w\in V$.

Define 
$Q_\Gamma:V\rightarrow \mathbb{F}_2$ by $${Q_\Gamma(v)}={g(v,v)}$$ for all $v\in V$.

\begin{proposition}\label{isquadratic}
$Q_\Gamma$ is a quadratic form with associated alternating form $f_\Gamma$, i.e., $Q_\Gamma$ is compatible with $\Gamma$.
\end{proposition}

%\begin{proof}
%We need to show that for all $v,w\in V$ we have
%$$Q_\Gamma(v\pluscirc w)=Q_\Gamma(v)+Q_\Gamma(w)+f_\Gamma(v,w).$$
%But this follows from the equalities
%$$\begin{array}{ll}
%(-1)^{Q_\Gamma(v\pluscirc w)}\one&=(v\pluscirc w)^2\\
%&=(vw)^2\\
%&=v(wv)w\\
%&=(-1)^{f_\Gamma(v,w)}v(vw)w\\
%&=(-1)^{f_\Gamma(v,w)}v^2w^2\\
%&=(-1)^{f_\Gamma(v,w)}\cdot (-1)^{Q_\Gamma(v)}\cdot (-1)^{Q_\Gamma(w)}\one\\
%&=(-1)^{f_\Gamma(v,w)+Q_\Gamma(v)+Q_\Gamma(w)}\one
%\end{array}$$

%\end{proof}

As the quadratic form $Q_\Gamma$ determines the alternating form $f_\Gamma$ as well as
the values of $\lambda$, we find the algebra $\mathfrak{C}(\Gamma)$ to be determined, up to isomorphism,  by the quadratic form $Q_\Gamma$.
For this reason we also denote the (isomorphism type of this) algebra by $\mathfrak{C}(V,Q_\Gamma)$. 
It is called a \emph{special} quasi-Clifford algebra. See also \cite{quasi}.

\section{The structure of the algebras}
\label{sect:structure}
\label{matrixalgebras}

As we have seen in the previous section, the algebras $\mathfrak{C}(\Gamma)$, where 
$\Gamma=(\mathcal{V},\mathcal{E},\lambda)$ is a labeled graph whose vertices are labeled by invertible elements 
from a commutative and  associative ring $R$ with $\pm 1$, are, up to isomorphism, algebras $\mathfrak{C}(V,g,\lambda)$ for some
bilinear form $g$ on an $\mathbb{F}_2$ vector space  $V=V_\Gamma$ and  map $\lambda$. 
Moreover the isomorphism type of the algebras is determined by the alternating form $f:V\times V\rightarrow \mathbb{F}_2$, defined by $f(v,w)=g(v,w)+g(w,v)$, and the map $\lambda$. If $\lambda$ takes only the values $\pm 1$, then 
the isomorphism type of the algebra is determined by a quadratic form $Q$ on $V$ with $f$ as associated bilinear form.

The classification of alternating and  quadratic forms on vector spaces of finite dimension over the field of $2$ elements is well known. We discuss this briefly.

Let $f$ be an alternating form on $V$, a finite dimensional vector space over $\mathbb{F}_2$.
The radical of $f$, defined as $\Rad(f)=\{v\in V\mid f(v,w)=0$ for all $w\in V\}$, is then a subspace of $V$,
and we can write $V$ as an orthogonal sum  $V=V_1\perp V_2\perp\cdots \perp V_k\perp  \Rad(f)$,
where $V_1,\dots V_k$ are $2$-dimensional hyperbolic spaces, i.e. they contain vectors $v,w$  with $f(v,w)=1$,
and $f(v_i,v_j)=0$ for $v_i\in V_i$ and $v_j\in V_j$ with $i\neq j$.
The form $f$ is called \emph{nondegenerate} if its radical is $\{0\}$.

If $Q$ is a quadratic form on $V$ with associated bilinear form $f$, then the radical of $f$ contains the radical of $Q$, defined as $\Rad(Q)=\{v\in\Rad(f)\mid Q(v)=0\}$, as a subspace of codimension at most $1$.
We call the form $Q$ \emph{nondegenerate} if and only if $\Rad(Q)=\{0\}$.

In dimension one there is, up to isomorphism, a unique nontrivial quadratic form  $Q$, with $Q(x)=x^2$ for $x\in \mathbb{F}_2$, which is  nondegenerate. It is called of $0$-type.
In dimension 2 we have, up to isomorphism, exactly two  nondegenerate forms $Q$. We have for $(x_1,x_2)\in \mathbb{F}_2^2$ (up to isomorphism)  $Q(x_1,x_2)=x_1x_2$, and $Q$ is called of $+$-type (also called hyperbolic type), or $Q(x_1,x_2)=x_1^2+x_1x_2+x_2^2$, and $Q$ called of $-$-type (or elliptic type).
In dimension $n>2$ we can distinguish, up to isomorphism, the following forms:

\hspace{0.5cm}\begin{minipage}{11cm}
\begin{enumerate}
 \item[$+$-type:] $V$ is an orthogonal sum $V_1\perp\cdots\perp V_k\perp \Rad(Q)$, where all $V_i$ are $2$-spaces of $+$-type. 
 \item[$-$-type:] $V$ is an orthogonal sum $V_1\perp\cdots\perp V_k\perp \Rad(Q)$, where all $V_i$ are $2$-spaces of $+$-type, except for one, which is of $-$-type. 
 \item[$0$-type:] $V$ is an orthogonal sum $V_1\perp\cdots\perp V_k\perp \Rad(Q)$, where all $V_i$ are $2$-spaces of $+$-type, except for one, which is one dimensional and of $0$-type. \\
 Notice, in this last case we find the radical of $f$ to be larger than the radical of $Q$.
\end{enumerate}
\end{minipage}

One of the key observations in the proof of this classification is that the type of a direct orthogonal sum of two nondegenerate spaces
is determined by the types of these summants. The orthogonal direct sum of  spaces of type $x$ and type $y$, where $x,y=\pm$,
gives us a space of type $x\cdot y$. 
We will frequently use these observations  in the sequel.
We note that the number of isomorphism classes of  quadratic spaces $(V,Q)$ over $\mathbb{F}_2$ of infinite dimension is much larger, see \cite{Hall_extra}.

The decomposition of $V$ into pairwise orthogonal subspaces provides a decomposition of the algebra $\mathfrak{C}(V,g,\lambda)$ into tensor products.
Indeed, if we assume $R$ to be a field, then the following proposition yields this decomposition.

\begin{proposition} \label{decompose}
Let $R$ be a field.
Suppose $V$ %is finite dimensional and 
can be decomposed as a direct orthogonal sum $V_1\perp V_2$.

For $i=1,2$, let $\mathcal{B}_i$ be a basis for $V_i$ and $g_i$ the restriction of $g$ to $V_i$.
Then $\mathfrak{C}(V,g,\lambda)$ is isomorphic to
$\mathfrak{C}(V_1,g_1,\lambda_1)\otimes \mathfrak{C}(V_2,g_2,\lambda_2)$, where for $i=1,2$ the map $\lambda_i:\mathcal{B}_i\rightarrow R^*$
is defined by $\lambda_i(b)=(-)^{g(b,b)}\Lambda(b,b)$ for all $b\in \mathcal{B}_i$.
\end{proposition}

\begin{proof}
The map $\phi$ that sends each tensor ${v_1}\otimes {v_2}\in \mathfrak{C}(V_1,g_1,\lambda_1)\otimes \mathfrak{C}(V_2,g_2,\lambda_2)$, with $v_1\in V_1, v_2\in V_2$ to ${v_1\pluscirc v_2}$ extends uniquely to a linear map $$\phi:\mathfrak{C}(V_1,g_1,\lambda_1)\otimes \mathfrak{C}(V_2,g_2,\lambda_2)\rightarrow \mathfrak{C}(V,g,\lambda).$$
Moreover, as the elements ${v_1}$ and ${v_2}$ commute in  $\mathfrak{C}(V,g,\lambda)$, it is straightforward to check that $\phi$ is a surjective homomorphism of algebras.
If $V_1$ and  $V_2$ are finite dimensional, then  the dimensions of $\mathfrak{C}(V,g,\lambda)$ and 
$\mathfrak{C}(V_1,g_1,\lambda_1)\otimes \mathfrak{C}(V_2,g_2,\lambda_2)$ coincide, and we find $\phi$ to be injective.

Now suppose $V$ is infinite dimensional.
Suppose $v\in \mathfrak{C}(V_1,g_1,\lambda_1)\otimes \mathfrak{C}(V_2,g_2,\lambda_2)$ is in the kernel of $\phi$.
Then as $v$ is a finite sum of scalar multiples of tensors of the form  $v_1\otimes v_2$, 
we can find finite dimensional subspaces $U_i\subseteq V_i$, for $i=1,2$, 
 which are spanned by  $\mathcal{B}_{U_i}=\mathcal{B}_i\cap{U_i}$. 
By the above, $\phi$ restricted to the tensor product $\mathfrak{C}(U_1,g_{1|U_1\times U_1},\lambda_{1|\mathcal{B}_{U_1}})\otimes \mathfrak{C}(U,g_{2|U_2\times U_2},\lambda_{2 |\mathcal{B}_{U_2}})$
is injective. So $v=0$.

This proves $\phi$ to be an isomorphism.
\end{proof}

In case $R=\mathbb{F}$ is a field, we can use the decomposition of the space $V$ into an orthogonal sum of $1$- and $2$-dimensional subspaces and  the above \cref{decompose} to find the structure of the algebra
$\mathfrak{C}(V,g,\lambda)$ for finite dimensional $V$.
The algebra $\mathfrak{C}(V,g,\lambda)$ is then  isomorphic to a tensor product of algebras $\mathcal{C}(V',g',\lambda')$, where $V'$ has dimension $1$ or $2$.
These $1$- and $2$-dimensional cases are well known and  worked out in for example \cite{quasi}.
For completeness we describe these algebras.

First assume $(V',g')$ is $1$-dimensional.
In this case $\mathfrak{C}(V',g',\lambda')$ is generated by a single element $v$, where $0\neq v\in V'$
satisfying $v^2=\lambda'(v)\one\in \mathbb{F}$.

If $\lambda'(v)$ is a square in $\mathbb{F}$, we find $\mathfrak{C}(V',g',\lambda')$ to be isomorphic to $\mathbb{F}\times \mathbb{F}$.
If $\lambda'(v)$ is not a square, then $\mathfrak{C}(V',g',\lambda')$ is isomorphic to the quadratic field extension $\mathbb{F}[\alpha]$, with $\alpha^2=\lambda'(v)$.

Now assume $V'$ to be a nondegenerate  $2$-dimensional space with respect to the alternating form $f'$.
So, $\mathfrak{C}(V',g',\lambda')$ is generated by two elements $v,w$ satisfying the relations
$$\begin{array}{ll}
v^2&=\alpha \one\\
w^2&=\beta \one\\
vw+wv&=0,\\
\end{array}$$
for some $\alpha,\beta\in \mathbb{F}^*$.

We are dealing with a quaternion algebra, which we denote by $\frac{(\alpha,\beta)}{\mathbb{F}}$ and can consider again two cases.

First, there exist $\gamma,\delta\in \mathbb{F}$ such that $\beta+\delta^2\alpha=\gamma^2$.
In this case we can identify  $\mathfrak{C}(V',g',\lambda')$ with the algebra of $2\times 2$ matrices over $\mathbb{F}$ via the map
$$v\mapsto \begin{pmatrix}0&\alpha\\1&0\end{pmatrix}\ \mathrm{and} \ w\mapsto \begin{pmatrix}\gamma&\alpha\delta\\-\delta&-\gamma\\\end{pmatrix}.$$

Second, there are no such $\gamma$ and $\delta$.
In this case the algebra is isomorphic to the subalgebra 
of the algebra of $4\times 4$ matrices over $\mathbb{F}$ generated by 
$$v\mapsto\begin{pmatrix} 0&\alpha &0&0\\1&0&0&0\\0&0&0&\alpha\\0&0&1&0\\\end{pmatrix}\ \mathrm{and}\ 
w\mapsto \begin{pmatrix} 0&0 &\beta&0\\0&0&0&-\beta\\1&0&0&0\\0&-1&0&0\\\end{pmatrix}.$$

This leads to the following result.

\begin{theorem}\label{tensorthm}
Let  $g$ be a bilinear form on 
a finite dimensional $\mathbb{F}_2$-space $V$ with basis $\mathcal{V}$, let $\mathbb{F}$ be a field,  and $\lambda:\mathcal{V}\rightarrow \mathbb{F}^*$ a map. 
Then $\mathfrak{C}(V,g,\lambda)$ over $\mathbb{F}$ is isomorphic to a tensor product of algebras
$\mathbb{F}\times \mathbb{F}$, quadratic field extensions of $\mathbb{F}$, the algebra of  $2\times 2$ matrices over $\mathbb{F}$, or
division algebras $\frac{(\alpha,\beta)}{\mathbb{F}}$, where $\alpha,\beta\in \mathbb{F}$ for which there are no
$\gamma$ and $\delta$ in $\mathbb{F}$ with $\beta+\delta^2\alpha=\gamma^2$.
\end{theorem}

\bigskip

In case $\lambda$ takes only the values $1$ and $-1$, we are dealing with a special quasi-Clifford algebra, and there is a quadratic form $Q$ with associated alternating  form $f$ such that
$\mathfrak{C}(\Gamma)$ is isomorphic to $\mathfrak{C}(V,Q)$, we can be more precise.
In this situation we  consider three types of fields, type I, II and III, defined by:

\hspace{0.5cm}\begin{minipage}{11cm}
\begin{enumerate}
\item[type I:] There is an element $i\in \mathbb{F}$ with $i^2=-1$.
  \item[type II:] There is no $i\in \mathbb{F}$ with $i^2=-1$, but there are $x,y\in\mathbb{F}$ with $x^2+y^2=-1$. 
\item[type III:] There are no $x,y\in \mathbb{F}$ with $x^2+y^2=-1$. 
  \end{enumerate}
\end{minipage}

If $V$ is $1$-dimensional, then
$\mathfrak{C}(V,Q)$ is isomorphic to $\mathbb{F}\times \mathbb{F}$ in case
$Q$ is trivial on $V$ or $\mathbb{F}$ is a field of type I.
If $Q$ is non-trivial on $V$ and $\mathbb{F}$ is of type
II or III, then $\mathfrak{C}(V,Q)$ is isomorphic to $\mathbb{F}[i]$,
where $i^2=-1$.

Now assume that $V=\langle e_1,e_2\rangle$ is $2$-dimensional and suppose
$Q$ is of $+$-type, $Q(e_1)=Q(e_2)=0$ and $f(e_1,e_2)=1$.
Then we can identify
$\mathfrak{C}(V,Q)$ with $M(2,\mathbb{F})$, the algebra of $2\times 2$-matrices
via the map $$e_1\mapsto \begin{pmatrix}0 & 1 \\ 1&0\end{pmatrix}$$
and $$e_2\mapsto \begin{pmatrix}1 & 0 \\ 0&-1\end{pmatrix}.$$

  If $Q$ is of $-$-type, then we may assume that $Q(e_1)=Q(e_2)=f(e_1,e_2)=1$  and
  we can identify $\mathfrak{C}(V,Q)$ with $M(2,\mathbb{F})$
via the map $$e_1\mapsto \begin{pmatrix}0&1\\-1&0\end{pmatrix}\textrm{ and }e_2\mapsto \begin{pmatrix}-i&0\\0&i\end{pmatrix}$$ if $\mathbb{F}$ is of type I,
     and 
     $$e_1\mapsto \begin{pmatrix}0&1\\-1&0\end{pmatrix}\textrm{ and }e_2\mapsto \begin{pmatrix}x&-y\\-y&-x\end{pmatrix}$$ if $\mathbb{F}$ is of type II and $x,y\in \mathbb{F}$ with $x^2+y^2=-1$.

         If $\mathbb{F}$ is of type III, we can identify $\mathfrak{C}(V,Q)$
         with the matrix algebra       
$$\left\langle \begin{pmatrix}0&-1&0&0\\1&0&0&0\\0&0&0&-1\\0&0&1&0\\ \end{pmatrix}, \begin{pmatrix}0&0&-1&0\\0&0&0&1\\1&0&0&0\\0&-1&0&0\\ \end{pmatrix} \right\rangle.$$ This algebra we denote by $\mathbb{H}$, the algebra 
    of quaternions over $\mathbb{F}$.

This implies that for finite dimensional orthogonal spaces $(V,Q)$
the algebra $\mathfrak{C}(V,Q)$ is determined, up to isomorphism, by the following parameters:

\begin{enumerate}[\rm(a)]
\item Dimension $n$ of $\overline{V}:=V/\mathrm{Rad}(Q)$;
\item Dimension $r$ of $\mathrm{Rad}(f)$;  
\item Type of $\overline{Q}$, the form induced by $Q$ on  $\overline{V}$;
\item Type of $\mathbb{F}$.
\end{enumerate}

We can now describe the various isomorphism classes of the algebras $\mathfrak{C}(V,Q)$ in terms of these parameters.

\begin{proposition}\label{martixalgebras}
  Let $(V,Q)$ be a nontrivial, finite dimensional quadratic space over the field $\mathbb{F}_2$.
  Then the isomorphism type of the special quasi-Clifford   algebra $\mathfrak{C}(V,Q)$ over a field $\mathbb{F}$ of characteristic $\neq 2$ is given in Table \ref{isotable}.
\end{proposition}

\medskip

\begin{table}[h]
\begin{tabular}{|l|l|l|l|}
  \hline
$\mathrm{dim}(\overline{V})$ &  Type($\overline{Q}$) &Type of $\mathbb{F}$& Algebra\\
\hline\hline
$n=0 \pmod{2}$&   $+$ & I& $(M(2,\mathbb{F})^{\otimes\frac{n}{2}})^{2^{r}}$ \\
$n=0 \pmod{2}$&   $-$ & I& $(M(2,\mathbb{F})^{\otimes\frac{n}{2}})^{2^{r}}$\\
$n=1 \pmod{2}$&  $0$   & I& $(M(2,\mathbb{F})^{\otimes\frac{n-1}{2}})^{2^{r}}$\\
$n=0 \pmod{2}$&   $+$ & II&  $(M(2,\mathbb{F})^{\otimes\frac{n}{2}})^{2^{r}}$ \\
$n=0 \pmod{2}$&   $-$ & II& $(M(2,\mathbb{F})^{\otimes\frac{n}{2}})^{2^{r}}$\\
$n=1 \pmod{2}$& $0$    & II& $(M(2,\mathbb{F})^{\otimes\frac{n-1}{2}}\otimes\mathbb{F}[i])^{2^{r-1}}$\\
$n=0 \pmod{2}$&   $+$ & III& $(M(2,\mathbb{F})^{\otimes\frac{n}{2}})^{2^{r}}$ \\
$n=0 \pmod{2}$&   $-$ & III& $(M(2,\mathbb{F})^{\otimes\frac{n-2}{2}}\otimes \mathbb{H})^{2^{r}}$\\
$n=1 \pmod{2}$& $0$    & III& $(M(2,\mathbb{F})^{\otimes\frac{n-1}{2}}\otimes\mathbb{F}[i])^{2^{r-1}}$\\
\hline
\end{tabular}

\medskip
\caption{The isomorphism types of the special quasi-Clifford algebras $\mathfrak{C}(V,Q)$, where $(V,Q)$ is a nontrivial quadratic space.}\label{isotable}
\end{table}

\section{Special quasi-Clifford algebras obtained from graphs of type $A$, $D$ and $E$}\label{sect:examples}
In this section we consider a few examples of algebras given by a graph $\Gamma$ whose vertices are labeled by $1$ and $-1$.
We only consider cases where the ring $R=\mathbb{F}$ is  a field of characteristic different from $2$. 
Vertices with a label $-1$ will be colored (and called) black and those with a value $1$ are colored (and called) white. 
When drawing a graph $\Gamma$ we use the  color gray for a  vertex to indicate that we have not yet determined whether its color should be black or white.

\begin{example}[Clifford algebras and graphs of type $A$]
Let $\Gamma$ be the complete graph on $n$ vertices with $p$ white vertices  and $q$ black vertices.
  Then of course $\mathfrak{C}(\Gamma)$ is isomorphic to the Clifford algebra $\mathrm{Cl}(p,q)$ over the field $\mathbb{F}$. 
  Consider the corresponding quadratic space $(V,Q)=(V_\Gamma,Q_\Gamma)$ obtained from $\Gamma$.  
  We embed $(V,Q)$ as a hyperplane in a  space $\widehat{V}$ with basis $e_0,\dots , e_n$ such that 
  the vectors $e_1\pluscirc e_0,\dots, e_n\pluscirc e_0$ correspond to the vertices of $\Gamma$,
  and extend both $Q$ and $f$ to $\widehat{V}$ such that $e_0$ is in the radical of   ${Q}$.
  So, we find $${Q}(e_i)=Q(e_i+e_0+e_0)=Q(e_i+e_0)+Q(e_0)=\begin{cases} 1 &\ \mathrm{if}\ 1\leq i\leq q\\
  0 &\ \mathrm{if}\ q+1\leq i\leq n\\
  \end{cases}.$$
  Moreover, $f(e_i,e_j)=1$ for all $i\neq j$ different from $0$.
 
 Clearly $V$ forms the subspace of $\widehat{V}$ of even weight vectors in the $e_i$.
 So, the vectors $f_i=e_i\oplus e_{i-1}$ with $1\leq i\leq n$ also form a basis  of $V$.
   Then $Q(f_i)=Q(e_i\pluscirc e_{i-1})=Q(e_i)+Q(e_{i-1})+1=1$ unless $i=q+1$, in which case it equals $0$.
  Moreover, $f(f_i,f_j)=1$ if and only if $|i-j|=1$.
  For this basis we find the graph $\Gamma'$ from \cref{newgraphAa} with $\mathfrak{C}(\Gamma)$ and $\mathfrak{C}(\Gamma')$ isomorphic. 
  (If $q=n$, then all vertices are black, if $q=0$, only $f_1$ is white.)

  This implies that for $q\geq 1$ we find  $\mathrm{Cl}(p,q)$ to be isomorphic to $\mathrm{Cl}(q-1,p+1)$.
  Just read the diagram from right to left.

  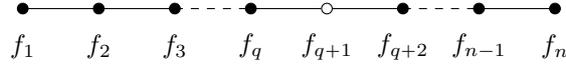
\begin{figure}[h]
    \begin{tikzpicture}

\filldraw[black] (0,0) circle [radius=2pt]
    (1,0) circle[radius=2pt]
    (2,0) circle[radius=2pt]
  (3,0) circle[radius=2pt]
  
  (5,0) circle[radius=2pt]
    (6,0) circle[radius=2pt]
    (7,0) circle[radius=2pt]
  ;

    \draw
    (0,-0.5) node {$f_1$}
    (1,-0.5) node {$f_2$}
    (2,-0.5) node {$f_3$}
    (3,-0.5) node {$f_{q}$}
    (4,-0.5) node {$f_{q+1}$}
    (5,-0.5) node {$f_{q+2}$}
    (6,-0.5) node {$f_{n-1}$}
    (7,-0.5) node {$f_n$};
    \draw (0,0)--(2,0);
    \draw (3,0)--(5,0);
    \draw (6,0)--(7,0);
    \draw[dashed] (2,0)--(3,0);
\draw[dashed] (5,0)--(6,0);

 \filldraw[white] (4,0) circle[radius=2pt];
\draw (4,0) circle[radius=2pt];

      \end{tikzpicture}
\caption{Graph of type $A_n$ obtained by changing the generators.}
\label{newgraphAa}
  \end{figure}

For $1\leq i\leq n$
let $$g_i=\begin{cases} e_i\pluscirc e_{i-1} \ \mathrm{for\ even} \ i\\
                        e_i\pluscirc\dots\pluscirc e_{0}\ \mathrm{for\ odd} \ i.\\
          \end{cases}  $$           

Then again we have a basis for $V$.
For even $i,j$ we have $f(g_i,g_j)=0$ and  $Q(g_i)=1$, unless $i=q+1$ in which case it is $0$.

For odd $i,j$ we also have $f(g_i,g_j)=0$.
Moreover, $Q(g_1)=1$ and for odd $i>1$  different from $q+1$  we find 
$$\begin{array}{ll}
Q(g_i)
&=Q(e_0\pluscirc \dots \pluscirc e_i)\\
&=Q(e_0\pluscirc \dots \pluscirc e_{i-2}\pluscirc e_{i-1}\pluscirc e_i)\\
&=Q(g_{i-2})+Q(e_{i-1}\pluscirc e_i)\\
&=Q(g_{i-2})+1.\end{array}
$$
If $i=q+1$, we find $Q(g_i)=Q(g_{i-2})$.

  From this information we can easily deduce the type of $Q$ restricted to each of the spaces $\langle g_i,g_{i+1}\rangle$,
  where $i$ is odd, and possibly $\langle g_{n}\rangle$ when $n$ is odd, and hence also the type of $Q$.
  In particular, we see that
  the type of $Q$ is multiplied with $-1$ if we raise
  $p$ or $q$ with $4$, and hence stays the same if we add $8$ to $q$ or $p$ (Bott-periodicity).
  Indeed, adding $4$ to $p$ or $q$ adds
%\begin{figure}[h]
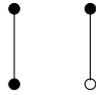
\begin{tikzpicture}[scale=0.5]
    \filldraw[black]

    (0,0) circle [radius=2pt]
    (0,1) circle [radius=2pt]
    (1,1) circle[radius=2pt];

    \draw
    (1,0) circle [radius=2pt];

    \draw
    (0,0)--(0,1)
    (1,0.05)--(1,1);
\end{tikzpicture}
%\end{figure}
  \noindent
  to the graph and multiplies the type of $Q$ with $-1$.

  For small values of $p$ and $q$
  we have collected this information in Table
  \ref{Cliffordtable}.

  \begin{table}
  \begin{tabular}{|l||l|l|l|l|l|}
    \hline
    $p\setminus q$ & 0  & 1 & 2   & 3 & 4 \\
    \hline\hline
    
         0         &  $+$ & 0 & $-$   &$-$  &  $-$  \\
         1         &  $+$ & $+$ & 0   &$-$  &  $-$   \\
         2         &  $+$ & $+$ & $+$   &$0$  & $-$  \\
         3         &  $0$ & $+$ & $+$   &$+$  &  $0$   \\
         4         &  $-$ & $0$ & $+$   &$+$  &  $+$  \\
    
    \hline
  \end{tabular}

  \medskip
  
  \caption{Type of $Q$ for small values of $p+q$.}
  \label{Cliffordtable}
  \end{table}

Using the results of Table \ref{Cliffordtable} and the above information, we find  in Table \ref{IsoCliffordTable}
 the isomorphism type of the Clifford algebras over fields $\mathbb{F}$ of type III.

  \begin{table}
    \begin{tabular}{|l|l|l|}
      \hline
      $p-q\pmod{8}$ & Type $Q$ & $\mathrm{Cl}(p,q)$\\
      \hline\hline
      0,2         & $+$        & $M(2^{\frac{p+q}{2}},\mathbb{F})$\\
        4,6         & $-$        & $M(2^{\frac{p+q-2}{2}},\mathbb{H})$\\

 1         & $+$        & $M(2^{\frac{p+q-1}{2}},\mathbb{F})^2$\\
 3,7         & $0$        & $M(2^{\frac{p+q-1}{2}},\mathbb{F}[i])$\\
      5         & $-$        & $M(2^{\frac{p+q-3}{2}},\mathbb{H})^2$\\     
        \hline

      \end{tabular}
    \medskip

    \caption{Isomorphism type of the Clifford algebras.}
    \label{IsoCliffordTable}
    \end{table}

  We notice that the above also classifies the algebras $\mathfrak{C}(\Gamma)$ where $\Gamma$ is a graph of type $A_n$,
  since we can replace the vertices $f_i$ by $e_i$, i.e., reserve the above described process, and end up with a complete graph. 
  
\end{example}

\begin{example}[Graphs of type $D$]
Next we consider graphs of type $D_n$, where $n\geq 4$. See Figure \ref{Graph_Dn}.
 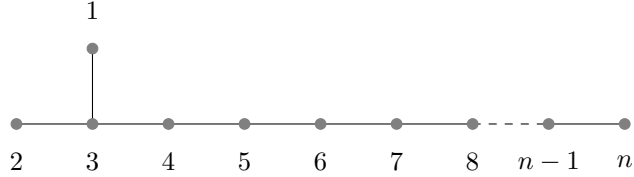
\begin{figure}[h!]
  \begin{tikzpicture}

    \draw
    (0,-0.5) node {$2$}
    (1,-0.5) node {$3$}
    (2,-0.5) node {$4$}
    (3,-0.5) node {$5$}
    (4,-0.5) node {$6$}
    (5,-0.5) node {$7$}
    (6,-0.5) node {$8$}
    (7,-0.5) node {$n-1$}
    (8,-0.5) node {$n$}
    (1,1.5) node {$1$};

    \draw (0,0)--(6,0);
    \draw (1,0)--(1,1);
    \draw (7,0)--(8,0);
    \draw[dashed] (6,0)--(7,0);
    \filldraw[gray]
    (0,0) circle[radius=2pt]
    (1,0) circle[radius=2pt]
    (2,0) circle[radius=2pt]
    (3,0) circle[radius=2pt]
    (4,0) circle[radius=2pt]
    (5,0) circle[radius=2pt]
    (6,0) circle[radius=2pt]
    (7,0) circle[radius=2pt]
    (8,0) circle[radius=2pt]
    (1,1) circle[radius=2pt];

  \end{tikzpicture}
  \caption{Graphs of type $D_n$.}
  \label{Graph_Dn}
  \end{figure}
  
 As we have seen in the above example, to classify the corresponding algebras we only have to consider the cases where at most one of the vertices $2,\dots, n$ is white.
 Moreover, we notice that $e_1\pluscirc e_2$ is an element which is in the radical of the form $f$ induced on $V=\langle e_1,\dots, e_n\rangle$.
 If both the vertices $1$ and $2$ are black or both are white, we find $e_1\pluscirc  e_2$ to be in the radical of $Q_\Gamma$ and $\mathfrak{C}(\Gamma)$ is the direct product  
 $\mathfrak{C}(\Gamma_1)\times \mathfrak{C}(\Gamma_1)$, where $\Gamma_1$ is obtained from $\Gamma$ by deleting vertex $1$.
 If only one of the two vertices $1$ and $2$ is black, then $Q_\Gamma(e_1+e_2)=1$ and we find $\mathfrak{C}(\Gamma)$ to be 
 isomorphic to  $\mathfrak{C}(\Gamma_1)\otimes \mathbb{F}[i]$.
\end{example}

\begin{example}[Graphs of type $E$]
  Let $\Gamma$ be a graph of type $E_n$, where $n\geq 1$ as in Figure \ref{En}.

  \begin{figure}[h!]
  \begin{tikzpicture}

    \filldraw[black]
    (0,0) circle[radius=2pt]
    (1,0) circle[radius=2pt]
    (2,0) circle[radius=2pt]
    (3,0) circle[radius=2pt]
    (4,0) circle[radius=2pt]
    (5,0) circle[radius=2pt]
    (6,0) circle[radius=2pt]
    (7,0) circle[radius=2pt]
    (8,0) circle[radius=2pt]
    (2,1) circle[radius=2pt];

    \draw
    (0,-0.5) node {$1$}
    (1,-0.5) node {$3$}
    (2,-0.5) node {$4$}
    (3,-0.5) node {$5$}
    (4,-0.5) node {$6$}
    (5,-0.5) node {$7$}
    (6,-0.5) node {$8$}
    (7,-0.5) node {$n-1$}
    (8,-0.5) node {$n$}
    (2,1.5) node {$2$};

    \draw (0,0)--(6,0);
    \draw (2,0)--(2,1);
    \draw (7,0)--(8,0);
    \draw[dashed] (6,0)--(7,0);
  \end{tikzpicture}
  \caption{Graphs of type $E_n$.}
  \label{En}
  \end{figure}
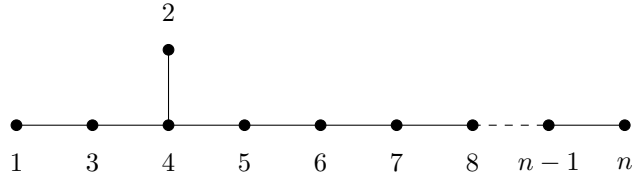

  Assume that all vertices are colored black.
  Consider the quadratic form on $V=\mathbb{F}_2^{n}$ given by
  $$Q(x_1,\dots,x_n)=(\sum_{i=1}^n\ x_i^2)+x_2x_4+x_1x_3+x_3x_4+\dots +x_{n-1}x_n.$$

  Then $Q(e_i)=1$ and $f(e_i,e_j)=1$ if and only if  $i$ is adjacent to $j$.
  So, $\mathfrak{C}(\Gamma)$ is isomorphic to $\mathfrak{C}(V,Q)$.

  For $n\geq 4$ even, we find that we can split $V$ into the orthogonal sum
  of the spaces $$\langle e_1,e_3\rangle\perp \langle e_2, e_2 \pluscirc e_4 \pluscirc \dots \pluscirc e_{n}\rangle \perp \langle e_5,e_6\rangle\perp \langle e_8, e_5\pluscirc  e_7\rangle\perp\dots\perp \langle e_n, e_5 \pluscirc e_7 \pluscirc\dots \pluscirc e_{n-1}\rangle.$$
 Such a $2$-dimensional space is of $+$ type if the second generator is of even weight, and of $-$ type if the second generator is of odd weight. 
  So we find $Q$ to be of $+$-type if $n=0,2\pmod{8}$ and of $-$-type for $n=4,6\pmod{8}$.

  For $n\geq 5$ odd  we find the vector $e_2\pluscirc e_5$ (for $n=5$) or $e_2 \pluscirc  e_5 \pluscirc e_7\pluscirc  e_9 \pluscirc\cdots \pluscirc e_n$ (for $n\geq 9$)
  to span the radical of $f$, the bilinear form associated to $Q$.
  This vector is isotropic if and only if $n=1\pmod{4}$.
  It remains to find the type of the form induced on $V/\mathrm{Rad}(f)$ in case $n=1\pmod{4}$.
  As modulo $e_2 \pluscirc e_5 \pluscirc e_7 \pluscirc e_9\pluscirc \cdots \pluscirc e_n$, we find that $e_2$ is in the subspace spanned by $e_1,e_3,e_4, \dots, e_n$, the
  type of $Q$ is determined by the type of $Q$ restricted to this subspace.
  As above we find that this is of $+$-type if $n-1=0,2\pmod{8}$ and of $-$-type if $n-1=4,6\pmod{8}$.
  So, also for graphs $E_n$ we find Bott-periodicity. 
  The information is summarized in \cref{EnCliffordTable}.

  \begin{table} 
  \begin{tabular}{|l|l|l|l|l|}
    \hline
    $\dim(V)$ & Type $Q$&Type of $\mathbb{F}$&$\mathfrak{C}(\Gamma)$\\
    \hline\hline
    $n=0\pmod{8}$ &$+$&I, II, III& $M(2,\mathbb{F})^{\otimes\frac{n}{2}}$\\
    $n=1\pmod{8}$ &$+$ &I, II, III& $(M(2,\mathbb{F})^{\otimes\frac{n-1}{2}})^2$\\
    $n=2\pmod{8}$ &$+$&I, II, III& $M(2,\mathbb{F})^{\otimes\frac{n}{2}}$\\
    $n=3\pmod{8}$ &$0$ &I& $(M(2,\mathbb{F})^{\otimes\frac{n-1}{2}})^2$\\
    $n=3\pmod{8}$ &$0$  &II,III& $M(2,\mathbb{F})^{\otimes\frac{n-1}{2}}\otimes\mathbb{F}[i]$\\
    $n=4\pmod{8}$ &$-$ &I,II& $M(2,\mathbb{F})^{\otimes\frac{n}{2}}$\\
    $n=4\pmod{8}$ &$-$ &III&$M(2,\mathbb{F})^{\otimes\frac{n-2}{2}}\otimes \mathbb{H}$ \\
    $n=5\pmod{8}$ &$-$ &I, II& $(M(2,\mathbb{F})^{\otimes\frac{n}{2}})^2$\\
    $n=5\pmod{8}$ &$-$ &III&$(M(2,\mathbb{F})^{\otimes\frac{n-3}{2}}\otimes \mathbb{H})^2$ \\
    $n=6\pmod{8}$ &$-$ &I, II& $M(2,\mathbb{F})^{\otimes\frac{n}{2}}$\\
    $n=6\pmod{8}$ &$-$ &III& $M(2,\mathbb{F})^{\otimes\frac{n-2}{2}}\otimes \mathbb{H}$\\
    $n=7\pmod{8}$ &$0 $ &I& $(M(2,\mathbb{F})^{\otimes\frac{n-1}{2}})^2$\\
    $n=7\pmod{8}$ &$0 $ &II,III& $M(2,\mathbb{F})^{\otimes\frac{n-1}{2}}\otimes\mathbb{F}[i]$\\
  \hline
  \end{tabular}

  \medskip
  
  \caption{Algebras $\mathfrak{C}(E_n)$.}  \label{EnCliffordTable}
\end{table}
  
\end{example}

\section{Involutions}
\label{sect:involutions}

Let $g$ be a bilinear form on an $\mathbb{F}_2$ vector space $V$ and denote by $f$ the alternating form on $V$
defined by $f(v,w)=g(v,w)+g(w,v)$ for all $v,w\in V$. 
Let $R$ be a commutative and associative ring containing distinct elements $1$ and $-1$.
Then consider the algebra $\mathfrak{C}:=\mathfrak{C}(V,g,\lambda)$ over $R$, where $\lambda$ is a map from a fixed basis $\mathcal{V}$ of $V$ to the set $R^*$
of invertible elements in the ring $R$.

In this section we describe some involutions acting on $\mathfrak{C}$, related to the grading, reversion and conjugation involutions of Clifford algebras.
We will describe these involutions with the help of linear and quadratic forms on $V$ with associated bilinear form $f$.

Let $H$ be a hyperplane of $V$ and define $\tau_H:\mathfrak{C}\rightarrow \mathfrak{C}$ by linear expansion of 
$$\tau_H(v)=\begin{cases}
v\ \mathrm{if}\ v\in H\\
-v\ \mathrm{if}\ v\not \in H.\\
\end{cases}$$

Let $Q$ be a quadratic form on $V$ with $f$ as associated  bilinear form.  
The involution $\tau_Q:\mathfrak{C}\rightarrow \mathfrak{C}$ is defined as the linear expansion of 

$$\tau_Q(v)=\begin{cases}
v\ \mathrm{if}\ Q(v)=0\\
-v\ \mathrm{if}\ Q(v)=1.\\
\end{cases}$$

\begin{proposition}\label{auto}
The involution $\tau_H$ is an automorphism of $\mathfrak{C}$.

The involution $\tau_Q$ is an anti-automorphism of $\mathfrak{C}$.
\end{proposition}

\begin{proof}
First consider $\tau_H$, where $H$ is a hyperplane of $V$.
It suffices to check for $u,v\in V\setminus \{0\}$ that $\tau_H(uv)=\tau_H(u)\tau_H(v)$.
As $H$ is a hyperplane, $\tau_H$ fixes either all three vectors $u,v,u\pluscirc v$ or negates two of them and, indeed, we find  $\tau_H(uv)=\tau_H(u)\tau_H(v)$.

To check that  $\tau_Q$ is an anti-automorphism, we have to check $\tau_Q(uv)=\tau_Q(v)\tau_Q(u)$.
If $Q(u)=Q(v)=0$, then $Q(u\pluscirc v)=0$ and $uv=vu$, or $Q(u\pluscirc v)=1$ and $uv=-vu$. In both cases $\tau_Q(uv)=\tau_Q(v)\tau_Q(u)$.

If $Q(u)=Q(v)=1$, then $Q(u\pluscirc v)=0$ and $uv=vu$ or $Q(u\pluscirc v)=1$ and $uv=-vu$. Again, in both cases $\tau_Q(uv)=\tau_Q(v)\tau_Q(u)$.

Finally, if $Q(u)=0$ and $Q(v)=1$ (or $Q(u)=1$ and $Q(v)=0$), then $Q(u\pluscirc v)=0$ and $uv=-vu$ or $Q(u\pluscirc v)=1$ and $uv=vu$. Also now we can check  $\tau_Q(uv)=\tau_Q(v)\tau_Q(u)$.
\end{proof}

\begin{proposition}\label{invol}
%Fix a quadratic form $Q$ with associated bilinear form $f$.
Let $\tau$ be a nontrivial linear map of $\mathfrak{C}$ mapping any $v\in V$ to $\pm v$.

If $\tau$ is an automorphism of $\mathfrak{C}$, then $\tau=\tau_H$ for some hyperplane $H$ of $V$.

If $\tau$ is an anti-automorphism of $\mathfrak{C}$, then $\tau=\tau_Q$ %or $\tau=\tau_Q\tau_H=\tau_H\tau_Q=\tau_{Q'}$ for some hyperplane $H$ of $V$ and quadratic form $Q'$ with associated bilinear form $f$.
for some quadratic form $Q$ with associated bilinear form $f$.
\end{proposition}

\begin{proof}
First assume that $\tau$ is an automorphism.
If $\tau$ fixes two vectors $v,w$, then it also fixes $v\pluscirc w$.
If $\tau$ negates two vectors $v,w\in V\setminus \{0\}$, then $v\pluscirc w$, should be fixed.
So, the vectors in $V$ fixed by $\tau$ form a hyperplane $H$ of $V$ and $\tau=\tau_H$.

Next, assume that $\tau$ is an anti-automorphism. Fix a quadratic from $Q$ with  associated bilinear form $Q$.
Then $\tau_Q\tau$ is an automorphism of $\mathfrak{C}$, and by the above, we either have $\tau=\tau_Q$ or $\tau_Q\tau_H=\tau_H\tau_Q$ for some hyperplane $H$ of $V$.

Suppose $\tau=\tau_Q\tau_H$ for some hyperplane $H$.
Let $h$ be the linear form on $V$ with kernel $H$ and define $Q'$ by $Q'(v)=Q(v)+h(v)$.
Then for all $v,w\in V$ we have $Q'(v\pluscirc w)=Q(v\pluscirc w)+h(v\pluscirc w)=Q(v)+Q(w)+h(v)+h(w)+f(v,w)=Q'(v)+Q'(w)+f(v,w)$.
So, $Q'$ is a quadratic form with associated bilinear form $f$.
Moreover, $\tau_{Q'}=\tau$. Indeed, if  $v\in H$, then $\tau(v)=\tau_Q\tau_H(v)=\tau_Q(v)=\tau_{Q+h}(v)$ and if $v\not \in H$, then
$\tau(v)=\tau_Q\tau_H(v)=-\tau_Q(v)=\tau_{Q+h}(v)$.
\end{proof}

Suppose $H$ is a hyperplane of $V$.
If $H$ contains the radical of $f$, then there is an element $v\in V$ with $H=\{w\in V\mid f(v,w)=0\}$.
But then $$\tau_v: \mathfrak{C}\rightarrow\mathfrak{C}$$ 
defined by  $\tau_v(w)= \frac{1}{\lambda(v)} vwv$  satisfies
$$\tau_v(w)=\frac{1}{\lambda(v)}vwv=\frac{1}{\lambda(v)} v^2w=w$$
for $w\in H$
and
$$\tau_v(w)=\frac{1}{\lambda(v)}vwv=-\frac{1}{\lambda(v)} v^2w=-w$$
for $w\not \in H$.
In particular,
we find  for all $w\in V$ that  
$$\tau_H(w)=\frac{1}{\lambda(v)}vwv=v^{-1}wv.$$

If $H$ does not contain the radical of $f$,
pick a point $r\not\in H$ with $r\in \Rad(f)$. 
Then $\mathfrak{C}$
can be written as a tensor product $\mathfrak{C}_{\langle r\rangle}\otimes \mathfrak{C}_{H}$,
where $\mathfrak{C}_{U}$ denotes the subalgebra of $\mathfrak{C}$ generated by $U\subseteq V$.

Now $\tau_H$ fixes  $\mathfrak{C}_{H}$ element wise, but acts nontrivially on  $\mathfrak{C}_{\langle r\rangle}$,
which is either isomorphic to $\mathbb{F}\times \mathbb{F}$ or to a quadratic extension 
$\mathbb{F}[\alpha]$ where $\alpha^2=\lambda(r)$.
In the first case, $\mathfrak{C}$ is the direct sum of two copies of $\mathfrak{C}_{H}$
and  $\tau$ fixes one copy and is multiplication with $-1$ on the other,
while in the second case  $\tau$ maps $\alpha$ to $-\alpha$ in $\mathbb{F}[\alpha]$,
so that $\tau$ acts as a field automorphism on $\mathfrak{C}$, viewed as an algebra over 
$\mathbb{F}[\alpha]$.

We also  analyse the action of the anti-automorphisms of \cref{invol} on $\mathfrak{C}$.

Let $v\in V$ with $Q(v)=1$. If $v^2=1$, then $v$ generates a subalgebra  $\mathbb{F}\cdot \one\times\mathbb{F}\cdot v$
on which $\tau_Q$ induces the map $$\lambda\cdot\one+\mu\cdot v\mapsto \lambda \cdot \one -\mu\cdot v.$$
Here $\lambda,\mu \in \mathbb{F}$.

Now let $v,w\in V$ with  $f(v,w)=1$.
Then $v,w$ generate  a subalgebra isomorphic to a quaternion algebra $\frac{(\alpha,\beta)}{\mathbb{F}}$.
If we define $Q$ in $V_i$ to be of $-$-type, with $Q(v)=Q(w)=1$, then
$\tau_Q$ induces quaternion conjugation on the subalgebra generated by $v$ and $w$:
$$\tau_Q(\lambda\cdot\one+\mu\cdot v+\rho\cdot w+\theta\cdot vw)= \lambda\cdot\one-\mu\cdot v-\rho\cdot w-\theta\cdot vw.$$
(Here $\lambda,\mu,\rho$ and $\theta$ are in $\mathbb{F}$.)

If $Q$ is of $+$-type, then up to permutation of $v,w$ and $vw$, we can assume that $\tau_Q$ 
only negates $vw$ and find  it to act as follows:
$$\tau_Q(\lambda\cdot\one+\mu\cdot v+\rho\cdot w+\theta\cdot vw)= \lambda\cdot\one+\mu\cdot v+\rho\cdot w-\theta\cdot vw.$$

If $r\in \Rad(f)$, and $Q(r)=1$, then, $\mathfrak{C}_{\langle r\rangle}$ is either
isomorphic to $\mathbb{F}\times\mathbb{F}$ and $\tau$ is multiplication by $1$ on one component and by $-1$ on the other,
or $\mathfrak{C}_{\langle r\rangle}$ is isomorphic to $\mathbb{F}[\alpha]$ where $\alpha^2=\lambda(r)$
and $\tau_Q$ induces the field automorphism mapping $\alpha$ to $-\alpha$.

In the finite dimensional case, we can write $V$ as an orthogonal sum
$$V=V_1\perp\dots \perp V_k\perp \Rad(Q),$$
as in described in \cref{sect:structure}, where $V_i$ are nondegenerate of dimension at most $2$.
Using this description and, conform \cref{decompose}, the observation that
$\mathfrak{C}$ is isomorphic to the tensor product of the corresponding algebras, we find $\tau_Q$
to be the tensor product of the various maps $\tau_{Q_{|V_i}}$ and the identity map on the subalgebra generated by $\Rad(Q)$.

\begin{example}\label{tau_on_matrix}
In the case that we are dealing with a special Clifford algebra, we can be even more precise, and give a concrete
description of the various anti-automorphisms.

We focus on the case that  $\mathfrak{C}$ is a finite dimensional special quasi-Clifford algebra $\mathfrak{C}(V,Q)$ where $\Rad(Q)=\{0\}$, defined over a field $\mathbb{F}$ of type III. In particular, we can identify $\mathfrak{C}$ with a matrix algebra over $\mathbb{F}$, over a quadratic extension $\mathbb{F}[i]$, or over a quaternion algebra
$\mathbb{H}$ as described in \cref{matrixalgebras}.

Then consider the anti-automorphism $\tau_Q$ acting on a non-degenerate $2$-space  of $+$-type, spanned by $v,w$ with $Q(v)=0=Q(w)$.
The subalgebra generated by $v,w$ can be identified with the algebra of $2\times 2$-matrices over $\mathbb{F}$ on which 
$\tau_Q$ induces matrix transposition.

So, if $Q$ is of $+$-type. Then $V$ is an orthogonal sum of $2$-space 
of dimension $2$ of $+$-type $\mathfrak{C}(V,Q)$ is a matrix algebra   on which $\tau_Q$  
induces matrix transposition.

If $Q$ is of $0$-type (or $-$-type), then $V$ is the orthogonal sum of a space of $+$-type and a $1$-dimensional space of $0$-type (or a $2$-space of $-$-type), and $\mathfrak{C}$ is 
a matrix algebra over a field $\mathbb{F}[i]$, where $i^2=-1$ (or  over the quaternions $\mathbb{H}=\frac{(-1,-1)}{\mathbb{F}}$).
We find $\tau_Q$ to induce matrix transposition followed by conjugation in $\mathbb{F}[i]$, (or quaternion conjugation in $\mathbb{H}$).

So, in all cases we find $\tau_Q$ induce the map
$$M\mapsto \overline{M}^\top$$
for $M$ in the matrix algebra over $\mathbb{F}$, $\mathbb{F}[i]$ and $\mathbb{H}$, respectively, where $\overline{\phantom{a}}$
is trivial or denotes conjugation in $\mathbb{F}[i]$ (or in $\mathbb{H}$).

If $Q'$ is a second quadratic form different from $Q$, then $\tau_{Q'}$ can be obtained as 
the composition $\tau_Q\tau_H=\tau_H\tau_Q$ for some hyperplane $H$ of $V$.

If $H$ contains $\Rad(f)$, where $f$ is the bilinear form associated to $Q$,  we find a $v\in V$ such that   
$$\tau_H(w)=v^{-1}wv$$
for  all $w\in V$.

We can have $Q(v)=0$ in which case we can assume $v$ to correspond to the matrix
$\begin{pmatrix}1&0\\0&-1\end{pmatrix}$
and find $\tau_{Q'}$ to induce the map
$$\begin{pmatrix}A&B\\C&D\\\end{pmatrix}\mapsto \begin{pmatrix}\overline{A}^\top&-\overline{C}^\top\\-\overline{B}^\top&\overline{D}^\top\\\end{pmatrix}.$$

If $Q(v)=1$, then we can assume  $v$ to correspond to the matrix
$\begin{pmatrix}0&1\\-1&0\end{pmatrix}$ and we find
$\tau_{Q'}$ to induce the map
$$\begin{pmatrix}A&B\\C&D\\\end{pmatrix}\mapsto\begin{pmatrix}\overline{D}^\top&-\overline{B}^\top\\-\overline{C}^\top&\overline{A}^\top\\\end{pmatrix}.$$

It remains to consider the case that $H$ does not meet $\Rad(f)$, which is non-trivial.
So, we are dealing with the case that $Q$ is of $0$-type.
In that case, as we have seen above, we find $\tau_H$ to be conjugation on $\mathbb{F}[i]$.
But then $\tau_{Q'}$ induces just matrix transposition on the matrix algebra over $\mathbb{F}[i]$. 

\end{example}
\section{Lie and Jordan algebras}
\label{sect:lie}

We continue with the notation of the previous sections.
Consider the algebra $\mathfrak{C}(V,g,\lambda)$ as in Section \ref{sect:twisted}, where $V$ is an $\mathbb{F}_2$ space equipped with a bilinear form $g$, with a fixed  basis $\mathcal{V}$ and $\lambda:\mathcal{V}\rightarrow R^*$, as  in 
Section \ref{sect:twisted}.  
This algebra can be identified with $\mathfrak{C}(\Gamma)$, where $\Gamma$ is the labeled graph with vertex set $\mathcal{V}$, two vertices $v,u$ adjacent
whenever $g(u,v)+g(v,u)=1$, and labeling $\lambda$. 

As
$\mathfrak{C}:=\mathfrak{C}(\Gamma)=\mathfrak{C}(V,g,\lambda)$ is an associative algebra, we can consider the associated Lie algebra $\mathfrak{C}^-$ and Jordan algebra $\mathfrak{C}^+$.

The Lie bracket $[\cdot,\cdot]$ is defined by the linear expansion of 

$$\begin{array}{ll}
  [v,w]&:=\frac{1}{2}(vw-wv)\\
  &=\begin{cases} vw& \mathrm{if}\ g(v,w)+g(w,v)=1,\\
    0&\mathrm{if}\ g(v,w)+g(w,v)=0,
    \end{cases}\\
    \end{array}$$
for all $v,w\in{V}$,
while the Jordan product $\circ$ is defined as

$$\begin{array}{ll}
v\circ w&:=\frac{1}{2}(vw+wv)\\
&=\begin{cases}
0 &\mathrm{if}\ g(v,w)+g(w,v)=1\\
 vw & \mathrm{if}\ g(v,w)+g(w,v)=0.\\
\end{cases}
\end{array}$$
for all $v,w\in{V}$.

Notice that  the Lie bracket $[\cdot,\cdot]$ as well as the Jordan product $\circ$ can be defined without assuming that $\frac{1}{2}\in R$.

Since the elements of $V$ form a basis for $\mathfrak{C}$, they also form a basis for $\mathfrak{C}^-$ and for $\mathfrak{C}^+$.

Clearly, the element $\one$ is in the center of the Lie algebra $\mathfrak{C}^-$,
but so are all elements $u\in {V}$ that are in the radical of $f$.
In particular we find that $\mathfrak{C}^-$ is a direct sum:
$$\mathfrak{C^-}= \mathfrak{r}\oplus \mathfrak{g}$$
where $\mathfrak{r}=\langle r\mid r\in \Rad(f) \rangle$ and $\mathfrak{g}=\langle v\mid v\in V\setminus \Rad(f)\rangle$
are Lie subalgebras of $\mathfrak{C}^-$.

The Lie algebra $\mathfrak{g}$ is also denoted by $\mathfrak{g}(V,g,\lambda)$ or $\mathfrak{g}(\Gamma)$.

\begin{proposition}\label{automorphism}
Let $H$ be a hyperplane of $V$, and $Q$ a quadratic form on $V$ with associated alternating form $f$. 
Then  $\tau_H$  and $-\tau_Q$ are  automorphisms of $\mathfrak{C}^-$ and $\mathfrak{g}$.

The maps $\tau_H$  and $\tau_Q$ are  automorphisms of $\mathfrak{C}^+$.
\end{proposition}

\begin{proof}
By \ref{auto} we find $\tau_H$  to be an automorphism of $\mathcal{C}$, and hence also of $\mathcal{C}^-$ and $\mathcal{C}^+$.
So, we consider $\tau_Q$, which is an anti-automorphism of $\mathcal{C}$.
Let $u,v\in V$. If $[u,v]=0$, then so is $[-\tau_Q(u),-\tau_Q(v)]$.
If $[u,v]\neq 0$, then

$$\begin{array}{ll}
-\tau_Q([u,v])&=-\tau_Q(uv)\\
&=-\tau_Q(v)\tau_Q(u)\\
&=-[\tau_Q(v),\tau_Q(u)]\\
&=[-\tau_Q(u),-\tau_Q(v)].
\end{array}$$

Similarly, if $u\circ v=0$, then so is $\tau_Q(u)\circ \tau_Q(v)$, and if $u\circ v\neq 0$, then 
% $u$ and $v$ commute and we find
$$\begin{array}{ll}
\tau_Q(u\circ v)&=\tau_Q(\frac{1}{2}(uv+vu))\\
&=\frac{1}{2}(\tau_Q(uv)+\tau_Q(vu))\\
&=\frac{1}{2}(\tau_Q(v)\tau_Q(u)+\tau_Q(u)\tau_Q(v))\\
&=\frac{1}{2}(\tau_Q(u)\tau_Q(v)+\tau_Q(v)\tau_Q(u))\\
&=\tau_Q(u)\circ\tau_Q(v).
\end{array}$$

\end{proof}

The centralizer in $\mathfrak{C}^-$ or $\mathfrak{g}$ (or $\mathfrak{C}^+$) of an automorphism $\sigma$ is of course 
a Lie subalgebra (Jordan algebra), which we denote by  $\mathfrak{C}^-_{\sigma}$ or $\mathfrak{g}_\sigma$, respectively (or $\mathfrak{C}^+_\sigma$).

Clearly, the centralizer $\mathfrak{g}_{\tau_H}$ of $\tau_H$, with $H$ a hyperplane of $V$, in $\mathfrak{g}=\mathfrak{g}(V,g,\lambda)$ is isomorphic to
$\mathfrak{g}(H,g_{\mid H},\lambda')$, where $\lambda'$ is the restriction of $\pm\Lambda$ to a basis of $H$. 
The subalgebra $\mathfrak{g}_{-\tau_Q}$  is the linear span of the set $\{ v\in V\mid Q(v)=1\}$ of non-isotropic vectors in $V$ inside $\mathfrak{g}$. 
Its subalgebra   $\mathfrak{g}_{Q}$  is defined to be the linear span of the set $\{ v\in V\setminus \Rad(f)\mid Q(v)=1\}$.

Similarly we find the centralizer of $\tau_H$ in $\mathfrak{C}^+$ to be isomorphic to the Jordan algebra
$\mathfrak{C}(H,g_{\mid H},\lambda')^+$ and 
the subalgebra $\mathfrak{C}^+_{\tau_Q}$  to be the linear span of the set $\{ v\in V\mid Q(v)=0\}$ of isotropic vectors in $V$ including $\one$.

\begin{lemma}\label{simple}
Suppose that $Q$ is non-degenerate.
If $V$ has dimension at least $3$, or $V$ has dimension $2$ and $Q$ is of $-$-type, 
then $\mathfrak{g}_{Q}$ 
is simple.

If  $V$ has dimension at least $3$, or $V$ has dimension $2$ and $Q$ is of $+$-type,
then the centralizer $\mathfrak{C}^+_{\tau_Q}$ of $\tau_Q$ in $\mathfrak{C}^+$ is simple.
\end{lemma}

\begin{proof}
Let $I$ be a nontrivial ideal in the centralizer in $\mathfrak{C}^-$ of $-\tau_Q$.
Suppose $0\neq r=\lambda_1x_1+\dots +\lambda_kx_k\in I$ with $x_i\in V$ with $\tau_Q(x_i)=-x_i$, and $x_i\not\in \Rad(Q)$ and $k$ minimal. 
Suppose $k>1$. 
Then take an element  $x\in V$ with $\tau_Q(x)=-x$ with 
$f(x, x_1)=0$ but  $f(x,x_2)=1$.
Then $[x,r]=\lambda_2[x, x_2]+\dots+\lambda_k [x,x_k]\in I$, which is a linear combination of elements from $V$ on which $Q$ takes the value $1$.
This contradicts the minimality of $k$.

Hence $k=1$. But then for each $x\in V$ with $Q(v)=f(x,x_1)=1$ we find $[x,[x,x_1]]$ and hence $x_1$ to be in $I$.
But then, when $V$ has dimension at least $3$, it is straightforward to check that all elements $v$ of $V$ with $Q(v)=1$ are in $I$, and $I$ is the full centralizer in $\mathfrak{C}^-$ of $-\tau_Q$.

The proof of the second statement is similar, but then focussing on vectors $v\in V$ with $Q(v)=0$. 
We only require that there are at least two isotropic $1$-spaces in $(V,Q)$, which is the case when 
$V$ has dimension at least $3$, or $Q$ is of $+$-type and $V$ has dimension $2$.

So, let $I$ be a nontrivial ideal in the centralizer in $\mathfrak{C}^+$ of $\tau_Q$.
Suppose $0\neq r=\lambda_1x_1+\dots +\lambda_kx_k\in I$ with $x_i\in V$ with $\tau_Q(x_i)=-x_i$ and $k$ minimal. 
Suppose $k>1$. 
Then take an element  $x\in V$ with $\tau_Q(x)=x$ with 
$f(x, x_1)=1$ but  $f(x,x_2)=0$.
Then $x\circ r=\lambda_2 x\circ x_2+\dots+\lambda_k x\circ x_k\in I$, which is a linear combination of elements from $V$ on which $Q$ takes the value $0$.
This contradicts the minimality of $k$.

Again we find $k=1$. But then $x_1^2$ and hence $\one\in I$ and $I=\mathfrak{C}^+_{\tau_Q}$.
So, $\mathfrak{C}^+_{\tau_Q}$ is simple.\end{proof}

\bigskip
\begin{example}
Consider a quasi-Clifford algebra  $\mathfrak{C}=\mathfrak{C}(V,g,\lambda)$
with symmetric bilinear form $f$  being  degenerate  on $V$. 

Let $r\in V$ be an element in the radical of $f$. Then $r$ commutes with all elements of $\mathfrak{C}$.
Morever, if $r^2=\one$ (which we can assume if $\Lambda(r,r)$ is a square) we can write $\mathfrak{C}$ as the direct sum $\mathfrak{I}_r^+\oplus\mathfrak{I}_r^-$,
where $\mathfrak{I}_r^\epsilon=\langle u+\epsilon ur\mid u\in V\rangle$ denotes the $\epsilon$-eigenspace, with $\epsilon=\pm 1$, of right (or left) multiplication with $r$.

For $u,w\in V$ we have
$$\begin{array}{ll}
[u+\epsilon ur,w]&=[u,w]+[\epsilon ur,w]\\
       &=[u,w]+\epsilon[u,w]r\\
\end{array}$$
and
$$\begin{array}{ll}
(u+\epsilon ur)\circ w&=u\circ w+\epsilon ur\circ w\\
            &=(u\circ w)+\epsilon(u\circ w)r\\
\end{array}$$
implying that $\mathfrak{I}_r^\epsilon$ is an ideal of both $\mathfrak{C}^-$ and $\mathfrak{C}^+$.
In particular, if $f$ is nontrivial,  $\mathfrak{I}_r^\epsilon$ is  proper ideal of both   $\mathfrak{C}^-$ and $\mathfrak{C}^+$.
\end{example}

\begin{example}
Assume that $Q$ is a non-degenerate form on $V$ and consider the special quasi-Clifford algebra $\mathfrak{C}=\mathfrak{C}(V,Q)$ over a field of $\mathbb{F}$ of type III. 
Then  we can identify $\mathfrak{C}(V,Q)$ with a matrix algebra over $\mathbb{F}$,
$\mathbb{F}[i]$, or $\mathbb{H}$.

By \cref{tau_on_matrix} we can distinguish   only a few  cases for an anti-automorphism 
which maps all elements  $v\in V$ to $\pm v$:
$\tau_1=\tau_Q$, $\tau_2=\tau_v\tau_Q$ with $v^2=\one$, $\tau_3=\tau_v\tau_Q$ with $v^2=-\one$, and
$\tau_4=\tau_Q\tau_H$ for some hyperplane not containing $\Rad(f)$, where $f$ is the bilinear form associated to $Q$.

This leads to the various isomorphism types of the corresponding Lie algebras $\mathfrak{g}_{-\tau}$ as provided in \cref{isoLietable}.

\begin{table}[h]
{\small

\begin{tabular}{|l|l|l|l|l|l|}
  \hline
$\mathrm{dim}(V)$ &  $Q$ &$\tau$ &$\mathfrak{C}$ &$\mathfrak{g}$&$\mathfrak{g}_{-\tau}$\\
\hline\hline
$n=0 \pmod{2}$&   $+$&$\tau_1$ &  $M(2,\mathbb{F})^{\otimes\frac{n}{2}}$& $\mathfrak{gl}(2^{\frac{n}{2}},\mathbb{F})$&$\mathfrak{so}(2^{\frac{n}{2}},\mathbb{F})$\\
$n=0 \pmod{2}$&   $+$&$\tau_2$ &  $M(2,\mathbb{F})^{\otimes\frac{n}{2}}$& $\mathfrak{gl}(2^{\frac{n}{2}},\mathbb{F})$&$\mathfrak{so}(2^{\frac{n-2}{2}},2^{\frac{n-2}{2}},\mathbb{F})$\\
$n=0 \pmod{2}$&   $+$&$\tau_3$ &  $M(2,\mathbb{F})^{\otimes\frac{n}{2}}$& $\mathfrak{gl}(2^{\frac{n}{2}},\mathbb{F})$&$\mathfrak{sp}(2^{\frac{n}{2}},\mathbb{F})$\\
$n=0 \pmod{2}$&   $-$&$\tau_1$ &  $M(2,\mathbb{F})^{\otimes\frac{n-2}{2}}\otimes \mathbb{H}$&$\mathfrak{gl}(2^{\frac{n-2}{2}},\mathbb{H})$&$\mathfrak{sp}(2^{\frac{n-2}{2}},\mathbb{H})$\\
$n=0 \pmod{2}$&   $-$&$\tau_2$ &  $M(2,\mathbb{F})^{\otimes\frac{n-2}{2}}\otimes \mathbb{H}$&$\mathfrak{gl}(2^{\frac{n-2}{2}},\mathbb{H})$&$\mathfrak{sp}(2^{\frac{n-4}{2}},2^{\frac{n-4}{2}},\mathbb{H})$\\
$n=0 \pmod{2}$&   $-$&$\tau_3$ &  $M(2,\mathbb{F})^{\otimes\frac{n-2}{2}}\otimes \mathbb{H}$&$\mathfrak{gl}(2^{\frac{n-2}{2}},\mathbb{H})$  &   $\mathfrak{so}(2^{\frac{n-2}{2}},\mathbb{H})$\\
$n=1 \pmod{2}$&   $0$ &$\tau_1$&  $M(2,\mathbb{F})^{\otimes\frac{n-1}{2}}\otimes\mathbb{F}[i]$&$\mathfrak{gl}(2^{\frac{n-1}{2}},\mathbb{F}[i])$&$\mathfrak{su}(2^{\frac{n-1}{2}},\mathbb{F}[i])$\\
$n=1 \pmod{2}$&   $0$ &$\tau_2$&  $M(2,\mathbb{F})^{\otimes\frac{n-1}{2}}\otimes\mathbb{F}[i]$&$\mathfrak{gl}(2^{\frac{n-1}{2}},\mathbb{F}[i])$&$\mathfrak{su}(2^{\frac{n-3}{2}},2^{\frac{n-3}{2}},\mathbb{F}[i])$\\
$n=1 \pmod{2}$&   $0$ &$\tau_3$&  $M(2,\mathbb{F})^{\otimes\frac{n-1}{2}}\otimes\mathbb{F}[i]$&$\mathfrak{gl}(2^{\frac{n-1}{2}},\mathbb{F}[i])$&$\mathfrak{so}(2^{\frac{n-1}{2}},\mathbb{F}[i])$\\
$n=1 \pmod{2}$&   $0$ &$\tau_4$&  $M(2,\mathbb{F})^{\otimes\frac{n-1}{2}}\otimes\mathbb{F}[i]$&$\mathfrak{gl}(2^{\frac{n-1}{2}},\mathbb{F}[i])$&$\mathfrak{so}(2^{\frac{n-1}{2}},\mathbb{F}[i])$\\
\hline
\end{tabular}
}
\medskip
\caption{The isomorphism types of the Lie algebras $\mathfrak{g}$ and $\mathfrak{g}_{-\tau}$. }\label{isoLietable}
\end{table}

\end{example}

\begin{remark}\label{Shirokov}
If we fix a hyperplane $H$ of $V$, then the group $\langle -\tau_Q,\tau_H\rangle$ is elementary abelian of order $2^2$.
The Lie algebra $\mathfrak{g}=\mathfrak{g}(V,g,\lambda)$ can be decomposed as
$$\mathfrak{g}=\mathfrak{g}_{1,1}\oplus \mathfrak{g}_{1,-1}\oplus\mathfrak{g}_{-1,1}\oplus\mathfrak{g}_{-1,-1},$$
where $\mathfrak{g}_{i,j}$ for $i,j=\pm 1$ denotes the intersection of the $i$-eigenspace of $-\tau_Q$ and $j$-eigenspace of $\tau_H$. 

Notice that for $i,j,k,l=\pm 1$ we have $$[\mathfrak{g}_{i,j},\mathfrak{g}_{k,l}]\subseteq \mathfrak{g}_{ik,jl}.$$

So, we find in $\mathfrak{g}$
Lie subalgebras $\mathfrak{g}_{1,1}\oplus \mathfrak{g}_{1,-1}$,  $\mathfrak{g}_{1,1}\oplus \mathfrak{g}_{-1,1}$ and  $\mathfrak{g}_{1,1}\oplus \mathfrak{g}_{-1,-1}$, which are just the centralizers of the involutions $-\tau_Q,\tau_H$ and $-\tau_Q\tau_H$  in $\langle -\tau_Q,\tau_H\rangle$.
Notice that $\tau_Q\tau_H=\tau_{Q'}$ where $Q'$ is the quadratic form defined by $Q'(v)=Q(v)+\phi_H(v)$ for all $v\in V$, with $\phi_H$ being the linear form on $V$ with kernel equal to $H$.
The form $Q'$ has also $f$ as its associated symplectic form.

These decompositions and the corresponding Lie subalgebras are investigated by Shirokov in \cite{Shirokov15,Shirokov16,Shirokov18} in case we are dealing with a real Clifford algebra.
 Actually, several results of \cite{Shirokov15,Shirokov16,Shirokov18} follow  directly from the above considerations and  \cref{isoLietable}.

 When  $\mathbb{F}$ is a field of type III, one can also consider the   $\mathbb{F}$-Lie subalgebras 
 $$\mathfrak{g}_{1,1}\oplus i\mathfrak{g}_{k,l}$$ (where $k,l=\pm 1$) of the  Lie algebra $\mathfrak{g}(V,Q)$ defined over $\mathbb{F}[i]$ with $i^2=-1$. 
 See also \cite{Shirokov15,Shirokov16,Shirokov18}.
 
 Of course, similar decompositions of the Jordan algebra $\mathfrak{C}^+$ can be considered.
 \end{remark}

\section{Lie algebras obtained from graphs}
\label{sect:liegraphs}

Let $\Gamma=(\mathcal{V},\mathcal{E},\lambda)$ be a labeled graph with all labels invertible elements in $R$. 
Then let $\mathfrak{C}^-(\Gamma)$ be the Lie algebra of $\mathfrak{C}(\Gamma)$ and by $\mathfrak{g}(\Gamma)$ its Lie subalgebra
generated by the vectors of $V_\Gamma$ outside the radical of $f_\Gamma$. (Here we use the notation of \cref{sect:relations}.)
The vertices in $\mathcal{V}$ do generate $\mathfrak{C}(\Gamma)$,
but need not generate the Lie algebra $\mathfrak{C}^-(\Gamma)$, nor the Lie subalgebra $\mathfrak{g}(\Gamma)$.

In this section we determine the Lie subalgebras of $\mathfrak{C}^-(\Gamma)$ generated by the vertices of the graph $\Gamma$.
We denote this subalgebra by $\mathfrak{K}(\Gamma)$.
If $\Gamma$ is not connected, then of course the Lie algebra 
$\mathfrak{K}(\Gamma)$ is the direct sum of the Lie algebras obtained from the connected components of $\Gamma$.
%Moreover, elements in the radical of $f_\Gamma$ are central in the Lie algebra  $\mathfrak{g}(\Gamma)$.
%So, for the remainder of this section we assume $\Gamma$ to be connected with a symplectic form $f_\Gamma$ that has a trivial radical.
So, we assume that $\Gamma$ is connected.

As in the previous sections we  identify $\mathfrak{C}:=\mathfrak{C}(\Gamma)$ with the algebra $\mathfrak{C}(V_\Gamma,g_\Gamma,\lambda)$. By $f_\Gamma$ we denote the bilinear form on $V_\Gamma$ defined by $f_\Gamma(v,w)=g_\Gamma(v,w)+g_\Gamma(w,v)$ for all $v,w\in V_\Gamma$. By $\mathfrak{g}$ we denote the Lie algebra $\mathfrak{g}(\Gamma)$.

Let $Q$ be the (unique) quadratic form on $V_\Gamma$ that takes the value $1$ on all vertices of $\Gamma$ and has $f_\Gamma$ as its associated alternating form.
Then $-\tau_Q$ fixes all vertices of $\Gamma$ and we find $\mathfrak{K}:=\mathfrak{K}(\Gamma)$ to be a Lie subalgebra 
of $\mathfrak{g}_{-\tau_Q}$, the subalgebra of $\mathfrak{g}$ centralized by $-\tau_Q$. 
So, we have the following inclusions of Lie algebras:
$$\mathfrak{K}\subseteq \mathfrak{g}_{-\tau_Q} \subseteq\mathfrak{g}\subseteq \mathfrak{C}^-.$$

The Lie subalgebra  $\mathfrak{K}$  generated by the vertices of $\Gamma$ will be studied with the help of the geometry of $(V_\Gamma,f_\Gamma)$.

Let $(V,f)=(V_\Gamma,f_\Gamma)$ be the symplectic space over $\mathbb{F}_2$ with addition $\pluscirc$.
If  $v,w\in V$ are two non-zero vectors with  $f(v,w)=0$, then the subspace $\langle v,w\rangle$ is called a \emph{singular line}.
If $v\neq w\in V$ are nonzero vectors with $f(v,w)=1$, then we call the $2$-dimensional subspace $\langle v,w\rangle$ a \emph{hyperbolic line} of $(V,f)$.
We identify these $2$-spaces with the set of three nonzero vectors $\{v,w,v\pluscirc w\}$ that they contain.
We will consider two partial linear spaces on  $V$.
The first one is the \emph{polar space} $\mathcal{S}(V,f)$, i.e. the partial linear space $(P,L)$,
where $P$ denotes the set of nonzero vectors of $V$ and $L$ is the set of singular lines (or triples of points in a singular line).
It satisfies the 'one or all' property: a point is collinear with one or all points of a line.
By $\mathcal{H}(V,f)$ we denote the partial linear space $(P',L')$ where $P'$ consists of all the vectors $v$ of $V\setminus \Rad(f)$
and whose lines  $L'$ are the hyperbolic lines. 
It is a so-called \emph{cotriangular space}, having the property that for each point $p$ and line $\ell$ not on $p$, the point $p$ is collinear to $0$ or all but one of the points of $\ell$.

A \emph{subspace} of a partial linear space $\Pi=(P,L)$ is a subset $S$ of the point set of $\Pi$  such that each line meeting $S$ in two points is contained in $S$. A subspace $S$ is often identified with the partial linear space $(S, \{\ell\in L\mid \ell\subseteq S\})$.
As the intersection of subspaces is again a subspace, we can define the subspace 
generated by a subset  $X$  of $P$ to be the intersection of all subspaces containing $X$.

Cotriangular spaces and their subspaces have been studied by several authors, see for example \cite{Hall,Shult, Seidel}.
Their connection with Lie algebras has been considered in \cite{cotriangleLie,Erik}.

\begin{proposition}
Suppose  ${\mathcal{S}}$ is a subspace of $\mathcal{H}(V_\Gamma,f_\Gamma)$.
Then $\mathcal{S}$ is a basis for the  Lie subalgebra of $\mathfrak{C}^-$ generated by  the elements in $\mathcal{S}$.
\end{proposition}
\begin{proof}
This follows immediately from the following observation: if $v,w\in V_\Gamma$ are collinear points in $\mathcal{S}$, then
$f_\Gamma(v,w)=1$. So $v\pluscirc w$ is a point of $\mathcal{S}$ and $[v,w]$ is a scalar multiple of $v\pluscirc w$.
If $v,w$ are not collinear, then $[v,w]=0$.
\end{proof}

The subspace of $\mathcal{H}(V,f)$ generated by $\mathcal V$ is  denote  by $\mathcal{H}(\Gamma)$. So, we find:

\begin{corollary}
The subspace ${\mathcal{H}}(\Gamma)$ of $\mathcal{H}(V_\Gamma,f_\Gamma)$ is  a basis for 
$\mathfrak{K}(\Gamma)$.
\end{corollary}

As we already noticed above, there is a unique quadratic form $Q$  that takes the value $1$ on all  $v\in \mathcal{V}$ and has $f$ as its associated
symplectic form. If $v,w\in V$ are two vectors with $Q(v)=Q(w)=1$ and $f(v,w)=1$, then also $Q(v\pluscirc w)=1$. So, the subset of all vectors $v\in V\setminus \Rad(f)$ with $Q(v)=1$ is a subspace  of  $\mathcal{H}(V,f)$, denoted by $\mathcal{H}(V,Q)$, containing $\mathcal{H}(\Gamma)$.

We have following inclusions:

$$\mathcal{H}(\Gamma)\subseteq \mathcal{H}(V,Q)\subseteq \mathcal{H}(V,f)\subseteq \mathbb{P}(V),$$
where $\mathbb{P}(V)$ denotes the projective space of $V$.

Below we will show that either $\mathcal{H}(\Gamma)$ equals  $\mathcal{H}(V,Q)$ or we have the following situation.

Let $\Omega$ and $\Omega'$ be two disjoint sets and denote by $P$ the set of all vectors
from the $\mathbb{F}_2$ vector space $W$ with basis $\Omega\cup \Omega'$ having exactly two $\Omega$-coordinates equal to $1$. 
By $L$ we denote the triples of $P$ having sum equal to $0$.
Then $(P,L)$ is a cotriangular space, whose points in $P$ generate $V_0$, the subspace of $W$ of all vectors with
an even number of $\Omega$-coordinates equal to $1$.
The standard dot product on $W$ induces the symplectic form $f_0$ on $V_0$ with radical spanned by
$\Omega'$ and, in case $\Omega$ is finite of even order the vector $r=\sum_{\omega\in \Omega}\ \omega$. 
Moreover, $Q_0$ is the unique quadratic form on $V_0$ with associated bilinear form $f_0$  
which takes the value $1$ on the elements in $P$.
Notice that, if $\Omega$ has even order, we find $Q_0(r)=1$.
So, the radical of $Q_0$ equals the span of $\Omega'$.

Now for any subspace $R_1$ of the radical $R_0$ of $Q_0$ which does not contain a vector $v+w$ with $v\neq w\in \Omega'$
the map  $p\in P\mapsto p+R_1$
provides an embedding of $(P,L)$ into $\mathcal{H}(V,Q)$,
where $V=V_0/R_1$ and $Q$ is the induced quadratic form by $Q_0$.

The cotriangular space $(P,L)$ is denoted by $\mathcal{T}(\Omega,\Omega')$ and its embedding into $\mathcal{H}(V,Q)$ is called a \emph{standard orthogonal embedding}.

The Lie subalgebra of $\mathfrak{C}^-$ generated by $P$ is then denoted by $\mathfrak{g}(\Omega,\Omega')$.

\begin{proposition}\label{subspacethm}
One of the following holds:
\begin{enumerate}
\item $\mathcal{H}(\Gamma)=\mathcal{H}(V,Q)$;
\item $\mathcal{H}(\Gamma)$ is a standard orthogonal embedding of  $\mathcal{T}(\Omega,\Omega')$ for some sets $\Omega$ and $\Omega'$ into $\mathcal{H}(V,Q)$.
\end{enumerate}
\end{proposition}
\begin{proof}
By \cite[Theorem 5.5]{Hall}
we have one of the following:

\begin{enumerate}
\item There is a symplectic form $f'$ on $V$ such that $\mathcal{H}(\Gamma)$ consists of all vectors outside the radical of $f'$;
\item There is a quadratic form $Q'$ on $V$ such that $\Pi$ consists of all vectors on which $Q'$ takes the value $1$.
\item ${\mathcal{H}(\Gamma)}$ is isomorphic to  $\mathcal{T}(\Omega,\Omega')$ for some sets $\Omega$, and $\Omega'$.
\end{enumerate}

In the first case we find for all basis elements $v,w\in \mathcal{V}$ that
$f(v,w)=f'(v,w)$ and we can conclude that $f=f'$.
But as $\mathcal{H}(\Gamma)$ is a proper subspace of $\mathcal{H}(V,f)$, this case can not occur.

Now consider the second case. 
Let $v,w\in \mathcal{H}(\Gamma)$ be collinear. Then $Q'(v)+Q'(w)+Q'(v\pluscirc w)=1=f(v,w)$.
If $v,w\in \mathcal{H}(\Gamma)$ are not collinear then $Q'(v)+Q'(w)+Q'(v\pluscirc w)=0=f(v,w)$.
So, for all  basis elements from $v,w\in \mathcal{V}$ we do have that  $Q'(v)+Q'(w)+Q'(v\pluscirc w)=f(v,w)$.
This clearly implies that $f$ is the alternating form associated to $Q'$.
But as $Q$ is the unique form which takes the value $1$ on all vertices of $\mathcal{V}$ and has $f$ as its associated
alternating form, we find $Q=Q'$.

Finally in the third case, where  ${\mathcal{H}(\Gamma)}$ is isomorphic to  $\mathcal{T}(\Omega,\Omega')$ for some sets $\Omega$, and $\Omega'$,
we find by using \cite[Theorem 5.4]{Hall} and \cite[Proposition 4.6]{Hallembed} that $V$ should be a quotient of the subspace $E$ of the $\mathbb{F}_2$-space with basis $\Omega\cup\Omega'$ of all vectors with
an even number of $\Omega$-coordinates equal to $0$, by some subspace $R$.
As in this quotient different points of $P$ should not be identified and the resulting space should carry the quadratic form $Q$, we find a standard embedding of  ${\mathcal{H}(\Gamma)}$.
\end{proof}

\begin{proposition}
$\mathcal{H}(\Gamma)$ is isomorphic to $\mathcal{H}(\Omega,\Omega')$
if and only if 
the graph $\overline{\Gamma}$ obtained  from $\Gamma$ by identifying vertices having the same set of neighbors, is  a line graph.
\end{proposition}

\begin{proof}
Assume  there are disjoint sets $\Omega, \Omega'$ such that
the vertices in $\Gamma$ can be identified with points of $\mathcal{T}(\Omega,\Omega')$.

As the vertices of $\Gamma$ span $V$, we find that for two vertices $v,w$ having the same set of neighbors in $\Gamma$
their sum is in the radical of $f$. Moreover, $Q(v+w)=Q(v)+Q(w)+f(v,w)=1+1+0=0$.
So, this sum is an element of $\mathbb{F}_2\Omega'$.

On the other hand if the sum of  vertices $v,w$ is in $\mathbb{F}_2\Omega'$ which is contained in the radical of $f$,
then these vertices will have the same set of neighbors in $\Gamma$.

This implies that after identification of the vertices with 
the same set of neighbors in $\Gamma$, we can  identify them with pairs of elements from $\Omega$, two such pairs being adjacent if and only if they intersect nontrivially, and find $\overline{\Gamma}$ to be a line graph.
\end{proof}

This implies the following.
\begin{theorem}\label{graphliealgebra}
The Lie subalgebra $\mathfrak{K}$ of $\mathfrak{g}$ is equal to $\mathfrak{g}_{-\tau_Q}$ if $\overline\Gamma$ is not a line graph 
or to $\mathfrak{g}(\Omega,\Omega')$ for some disjoint sets $\Omega,\Omega'$ with a standard orthogonal embedding of 
$\mathcal{T}(\Omega,\Omega')$ into $\mathcal{H}(V,Q)$ if $\overline\Gamma$ is a line graph.
\end{theorem}

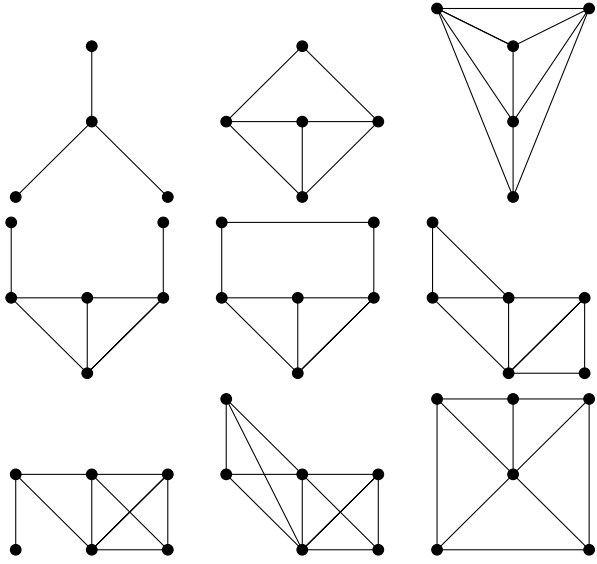
\begin{figure}
 \begin{tikzpicture}

    \filldraw[black]
    (0,0) circle[radius=2pt]
    (0,1) circle[radius=2pt]
    (-1,-1) circle[radius=2pt]
    (1,-1) circle[radius=2pt];

    \draw (0,0)--(0,1);
    \draw (-1,-1)--(0,0)--(1,-1);

  \end{tikzpicture}\hspace{0.5cm}
\begin{tikzpicture}

    \filldraw[black]
    (0,0) circle[radius=2pt]
    (0,1) circle[radius=2pt]
    (-1,0) circle[radius=2pt]
    (0,-1) circle[radius=2pt]
    (1,0) circle[radius=2pt];
    
    \draw (-1,0)--(1,0)--(0,1)--(-1,0)--(0,-1)--(1,0);
    \draw (0,0)--(0,-1);

  \end{tikzpicture}\hspace{0.5cm}
 \begin{tikzpicture}

    \filldraw[black]
    (0,1) circle[radius=2pt]
    (0,0) circle[radius=2pt]
    (0,-1) circle[radius=2pt]
    (1,1.5) circle[radius=2pt]
    (-1,1.5) circle[radius=2pt];
    
    \draw (-1,1.5)--(1,1.5)--(0,-1)--(-1,1.5)--(0,1)--(0,-1);
    \draw (-1,1.5)--(0,0)--(1,1.5);
\draw (-1,1.5)--(0,1)--(1,1.5);

  \end{tikzpicture}
   
  \medskip

\begin{tikzpicture}

    \filldraw[black]
    (0,0) circle[radius=2pt]
    (1,1) circle[radius=2pt] 
    (-1,1) circle[radius=2pt]
    (-1,0) circle[radius=2pt]
    (0,-1) circle[radius=2pt]
    (1,0) circle[radius=2pt];
    
    \draw (-1,1)--(-1,0)--(1,0)--(0,-1)--(1,0)--(1,1);
    \draw (0,0)--(0,-1)--(-1,0);

  \end{tikzpicture}\hspace{0.5cm}
\begin{tikzpicture}

    \filldraw[black]
    (0,0) circle[radius=2pt]
    (1,1) circle[radius=2pt] 
    (-1,1) circle[radius=2pt]
    (-1,0) circle[radius=2pt]
    (0,-1) circle[radius=2pt]
    (1,0) circle[radius=2pt];
    
    \draw (-1,1)--(-1,0)--(1,0)--(0,-1)--(1,0)--(1,1);
    \draw (0,0)--(0,-1)--(-1,0);
    \draw (-1,1)--(1,1);
  \end{tikzpicture}\hspace{0.5cm}
\begin{tikzpicture}

    \filldraw[black]
    (0,0) circle[radius=2pt]
    (1,-1) circle[radius=2pt] 
    (-1,1) circle[radius=2pt]
    (-1,0) circle[radius=2pt]
    (0,-1) circle[radius=2pt]
    (1,0) circle[radius=2pt];
    
    \draw (-1,1)--(-1,0)--(1,0)--(0,-1)--(1,0)--(1,-1);
    \draw (0,0)--(0,-1)--(-1,0);
    \draw (-1,1)--(0,0);
    \draw (0,-1)--(1,-1);
  \end{tikzpicture}
\medskip

  \begin{tikzpicture}

    \filldraw[black]
    (0,0) circle[radius=2pt]
    (1,-1) circle[radius=2pt] 
    (-1,-1) circle[radius=2pt]
    (-1,0) circle[radius=2pt]
    (0,-1) circle[radius=2pt]
    (1,0) circle[radius=2pt];
    
    \draw (-1,0)--(1,0)--(0,-1)--(1,0)--(1,-1);
    \draw (0,0)--(0,-1)--(-1,0);
   
    \draw (0,-1)--(1,-1);
    \draw (-1,-1)--(-1,0);
    \draw (1,-1)--(0,0);
  \end{tikzpicture}\hspace{0.5cm}
\begin{tikzpicture}

    \filldraw[black]
    (0,0) circle[radius=2pt]
    (1,-1) circle[radius=2pt] 
    (-1,1) circle[radius=2pt]
    (-1,0) circle[radius=2pt]
    (0,-1) circle[radius=2pt]
    (1,0) circle[radius=2pt];
    
    \draw (-1,1)--(-1,0)--(1,0)--(0,-1)--(1,0)--(1,-1);
    \draw (0,0)--(0,-1)--(-1,0);
    \draw (-1,1)--(0,0);
    \draw (0,-1)--(1,-1);
    \draw (-1,1)--(0,-1);
    \draw (1,-1)--(0,0);
  \end{tikzpicture}\hspace{0.5cm}
  \begin{tikzpicture}

    \filldraw[black]
    (0,0) circle[radius=2pt]
    (1,1) circle[radius=2pt] 
    (0,1) circle[radius=2pt]
    (-1,1) circle[radius=2pt]
    (-1,-1) circle[radius=2pt]
    (1,-1) circle[radius=2pt];
    
    \draw (1,1)--(0,1)--(-1,1)--(-1,-1)--(1,-1)--(1,1);
    \draw (0,0)--(1,1);
    \draw (0,0)--(0,1);
    \draw (0,0)--(-1,1);
    \draw (0,0)--(-1,-1);
    \draw (0,0)--(1,-1);
  \end{tikzpicture}

\caption{The nine forbidden subgraphs for a line graph.}\label{forbidden}
\end{figure}

We notice that due to Beineke's characterization of line graphs, see \cite{Beineke}, we can conclude that
$\overline{\Gamma}$ is not a line graph if it contains an induced subgraph $\Delta$ which is one of the  nine graphs from Figure \ref{forbidden}.
The three graphs on the first row of Figure \ref{forbidden}
are not reduced, while the others are.
(Here a graph $\Delta$ is called reduced if  $\overline{\Delta}= \Delta$.)
So, if $\Delta$ is one of these three graphs contained as an induced subgraph of $\overline{\Gamma}$, then 
$\overline\Gamma$ contains a vertex distinguishing the vertices that have in $\Delta$ 
the same set of neighbors.
So, in $\overline\Gamma$ we find  two vertices if $\Delta$ is  the first graph and one vertex in case $\Delta$ is the second or third  graph, such that adding these vertices to $\Delta$ we obtain a reduced graph.

This implies that $\Gamma$ contains a reduced connected subgraph $\Gamma_0$ on $6$ vertices which is not a line graph.
In particular, if we determine  the quadratic space $(V_{\Gamma_0},Q_{\Gamma_0})$ for this subgraph,
then this is an orthogonal $\mathbb{F}_2$-space of dimension  $6$.
If it is nondegenerate and of $+$-type, then its cotriangular space $\mathcal{H}(V_{\Gamma_0},Q_{\Gamma_0})$ is isomorphic to
$\mathcal{T}(\Omega)$, where $\Omega$ is of size $8$, contradicting that $\Gamma_0$ is not a line graph.
If it is degenerate with radical $R$, then we can view $\Gamma_0$ to be embedded in the cotriangular space of the 
quotient space $V_{\Gamma_0}/R$. Indeed, as $\Gamma_0$ is reduced, no two vertices are equal modulo $R$. The  cotriangular space of $V_{\Gamma_0}/R$ is isomorphic to $\mathcal{T}(\Omega)$ for some set $\Omega$ of size at most $6$, and we find $\Gamma_0$ to be a line graph as well, again leading to a contradiction.
So, $\Gamma_0$ generates a nondegenerate $6$-dimensional orthogonal space of $-$-type.
We have proven the following.

\begin{proposition}\label{6graph}
Suppose $\Gamma$ is a connected graph such that $\overline{\Gamma}$ is not a line graph.
Then $\Gamma$ contains a subgraph $\Gamma_0$ on $6$ vertices spanning a nondegenerate $6$-dimensional orthogonal $\mathbb{F}_2$ space  $(V_{\Gamma_0},Q_{\Gamma_0})$ of $-$-type.
\end{proposition}

\begin{corollary}\label{6graphcor}
Let $\Gamma=(\mathcal{V},\mathcal{E},\lambda)$ be a connected labeled graph with labels $\pm 1$ in a field $\mathbb{F}$ of odd characteristic.
Assume that   $\overline{\Gamma}$ is not a line graph.
Then, there is a subset $\mathcal{V}_0$ of size $6$ in $\mathcal{V}$ such that
the Lie subalgebra of $\mathfrak{g}(\Gamma)$ generated by $\mathcal{V}_0$ is isomorphic to
$\mathfrak{sp}(8,\mathbb{F})$ if $\mathbb{F}$ is of type I,
$\mathfrak{so}(8,\mathbb{F})$ if $\mathbb{F}$ is of type II,
or  $\mathfrak{sp}(4,\mathbb{H})$,
where $\mathbb{H}$ is a quaternion algebra over $\mathbb{F}$, if $\mathbb{F}$ is of type III.
\end{corollary}

\begin{proof}
Let $\Gamma_0$ be the subgraph on $6$ vertices guaranteed by
\cref{6graph}. Then $\mathfrak{K}(\Gamma_0)$ is the subalgebra we are looking for.
\end{proof}

\begin{remark}
The  Lie algebra $\mathfrak{sp}(4,\mathbb{H})$, where $\mathbb{H}$ is the division ring of real quaternions,
is the maximal compact Lie subalgebra of a split real Lie algebra of type $\mathfrak{e}_6$.
See \cref{compact}.
\end{remark}

\begin{remark}
Combining a polynomial algorithm to compute $\overline{\Gamma}$ from a graph $\Gamma$ (see for example Algorithm 2 in \cite{partition}),
and a polynomial algorithm  for checking whether a graph 
$\Gamma$ is a line graph or not (see for example \cite{alglinegraph}), 
yields a polynomial algorithm to determine whether $\mathfrak{K}(\Gamma)$ equals $\mathfrak{g}_{-\tau_Q}$
or $\mathfrak{g}(\Omega,\Omega')$ for some sets $\Omega$ and $\Omega'$.
\end{remark}

\begin{remark}
\cref{6graph} and \cref{6graphcor} are closely related to some results of Seven \cite{seven}.
See in particular \cite[Theorem 2.7]{seven}. Seven shows, among other things, the following:

Let $(V,Q)$ be an orthogonal space over the field with two elements with corresponding bilinear
form $f$. To each vector $v$ with $Q(v)=1$ we can assign a transvection $\tau_v:V\rightarrow V$
in the orthogonal group of $\mathrm{O}(V,Q)$, such that for all $w\in V$ we have $$\tau_v(w)=w+f(v,w)v.$$
Let $\mathcal{V}$ be a basis of anisotropic vectors of $V$, and denote by $\Gamma$ the graph where
two elements $v,w\in\mathcal{V}$ are adjacent if and only if $f(v,w)=1$.
If $\Gamma$ is connected, but $\overline{\Gamma}$ is not a line graph,
then $\Gamma$ contains an induced subgraph $\Gamma_0$ on six points that generate
a  nondegenerate $6$-dimensional orthogonal $\mathbb{F}_2$ space  $(V_{\Gamma_0},Q_{\Gamma_0})$ of $-$-type on which the corresponding transvections induce
the  orthogonal group $\mathrm{O}(V_{\Gamma_0},Q_{\Gamma_0})$ which is isomorphic to the Weyl group of type $E_6$.
\end{remark}

\begin{remark}
The Lie algebras considered in this and the previous section can be defined over $\mathbb{Z}$ and have a $\mathbb{Z}$-basis.
This implies that we can also consider them over $\mathbb{F}_2=\mathbb{Z}/2\mathbb{Z}$.
In this case we encounter the Kaplansky Lie algebras as described in \cite{Kaplansky}.
The involutions $\tau_Q$ defined over $\mathbb{Z}$ are closely related to the gradings described in
\cite{bouarroudj}.
\end{remark}

\section{A characterization of Lie algebras associated to a graph}\label{sect:Liegraph}

Suppose $\Gamma=(\mathcal{V},\mathcal{E},\lambda)$ is a labeled graph as in the previous sections.
Then the generators $x\neq y\in \mathcal{V}$ of the Lie algebra $\mathfrak{K}:=\mathfrak{K}(\Gamma)$ do satisfy the relations

$$\begin{array}{rll}
[x,y]&=0 & \textrm{ if }(x,y)\not\in \mathcal{E}\\
{[x,[x,y]]}&=\lambda(x) y & \textrm{ if }(x,y)\in \mathcal{E}.\\
\end{array}$$

So, the free Lie algebra $\mathfrak{g}_\Gamma$ with generators in $\mathcal{V}$ subject to the above relations
has then $\mathfrak{K}$ as a quotient.

\begin{theorem}\label{universallie}
Let $\Gamma=(\mathcal{V},\mathcal{E},\lambda)$ be a labeled graph with labels from $R^*$ containing the element $\frac{1}{2}$.

Suppose $\mathfrak{A}$ is an associative algebra with identity $\one_{\mathfrak{A}}$
and $\mathfrak{g}_0$  is a subalgebra of $\mathfrak{A}^-$
generated by elements $\overline v$ with $v\in \mathcal{V}$ 
with $$\begin{array}{lll}
\overline{v}^2&=\lambda(v)\one_{\mathfrak{A}}&\\
{[\overline{v},\overline{w}]}&=0&\mathrm{if}\ v\not \sim w\\
{[\overline{v},[\overline{v},\overline{w}]]}&=\lambda(v)\overline{w}&\mathrm{if}\ v\sim w.\\
\end{array}$$
Then $\mathfrak{g}_0$ is isomorphic to a quotient of  $\mathfrak{K}(\Gamma)$.
\end{theorem}

\begin{proof}
  For $v,w\in \mathcal{V}$ we have
  $ \overline{v}\overline{w}-\overline{w}\overline{v}=0$ or $\frac{1}{4}(\overline{v}(\overline{v}\overline{w}-\overline{w}\overline{v})-(\overline{v}\overline{w}-\overline{w}\overline{v})\overline{v})=\lambda(v) \overline{w}$.
  Suppose we are in the latter case.
  Then
  $\overline{v}^2\overline{w}-2\overline{v}\overline{w}\overline{v}+\overline{w}\overline{v}^2=4\lambda(v)\overline{w}$ and hence
  $2\lambda(v) \overline{w}+2\overline{v}\overline{w}\overline{v}=0$.
  Now multiplying with $\overline{v}$ yields
  $2\lambda(v)(\overline{v}\overline{w}+\overline{w}\overline{v})=0$ from which we deduce $\overline{v}\overline{w}+\overline{w}\overline{v}=0$.
  
  So, the  elements $\overline{v}$ with $v\in \mathcal{V}$ satisfy, as elements from $\mathfrak{A}$, the relations of the generators in $\mathcal{V}$ of $\mathfrak{C}(\Gamma)$.
  But then the subalgebra of $\mathfrak{A}$ generated by these elements is isomorphic to a quotient of $\mathfrak{C}(\Gamma)$.
  This implies that $\mathfrak{g}_0$ is isomorphic to a quotient of $\mathfrak{K}(\Gamma)$.
  \end{proof}

\cref{universallie} together with \cref{graphliealgebra} imply \cref{mainlie}.

The above result is motivated by, and a generalisation of, the results of \cite{spin}.
These connections are discussed in the next section.

\section{Spin representations and compact subalgebras of Kac-Moody algebras}
\label{sect:kac}

Let $\Gamma=(\mathcal{V},\mathcal{E},\lambda)$ be a labeled graph.
In \cite{spin}  linear representations $\phi$ of the free Lie algebra $\mathfrak{g}_\Gamma$  into $\mathfrak{gl}(W)$, the general linear Lie algebra
on a vector space $W$ over a field of characteristic not $2$, are considered.
If $x,y$ are two linear maps on $W$, then by $xy$ we denote the composition, and we consider the Lie product of $\mathfrak{gl}(W)$ to be defined
as $$[x,y]=\frac{1}{2}(xy-yx).$$

Such a representation $\phi$ is called a \emph{$\Gamma$-spin representation} of $\mathfrak{g}_\Gamma$, if and only if
$$\phi(x)^2=\lambda(x)\one_W$$
for all generators $x\in \mathcal{V}$.
Our first observation is that $\mathfrak{g}_\Gamma$ always admits such a representation.

\begin{proposition}\label{spin_exists}
The Lie algebras $\mathfrak{g}_\Gamma$ and  $\mathfrak{K}(\Gamma)$ admit a $\Gamma$-spin representation.
\end{proposition}

\begin{proof}
As $\mathfrak{K}(\Gamma)$ is a quotient of $\mathfrak{g}_\Gamma$, we only have to show that $\mathfrak{K}(\Gamma)$ admits a $\Gamma$-spin representation.

As the elements of $\mathfrak{K}(\Gamma)$ act by left multiplication on $\mathfrak{C}(\Gamma)$, and $v^2=\lambda(v)\one$ for all $v\in \mathcal{V}$,
we have found a $\Gamma$-spin representation of $\mathfrak{K}(\Gamma)$.
\end{proof}

 For finite graphs $\Gamma$ in which  all vertices are labeled by $-1$, a $\Gamma$-spin representation is called a \emph{generalized spin representation}
 in \cite{spin,Buyl,Damour}.
 Such representations are constructed in \cite{spin}, generalizing \cite{Buyl, Damour}
 in which such representations have been constructed for graphs of type $E_{9}$ and $E_{10}$.

The following characterization of the $\Gamma$-spin representations generalizes the work in \cite{spin}.

\begin{theorem}\label{spinthm}
  Suppose $\phi:\mathfrak{g}_\Gamma\rightarrow \mathfrak{gl}(W)$ for some vector space $W$ over  a field of characteristic $\neq 2$ is a 
  $\Gamma$-spin representation of $\mathfrak{g}_\Gamma$.
  
  Then $\phi(\mathfrak{g}_\Gamma)$ is isomorphic to a quotient of $\mathfrak{K}(\Gamma)$.
\end{theorem}

\begin{proof}
This follows directly by \cref{universallie}.
\end{proof}

  \bigskip

 Over fields $\mathbb{F}$ of characteristic $0$ Berman \cite{Berman} relates the Lie algebra $\mathfrak{g}_\Gamma$, for graphs $\Gamma$ with all vertices  labeled with $-1$, to the so-called compact subalgebras of Kac-Moody algebras.
Let us explain this connection, restricting ourselves to the simply laced case.

Let $A=(a_{ij})$ be a generalized Cartan matrix indexed by the set $\mathcal{V}$, which is simply laced. That means  $$a_{ii}=2$$ and 
$$a_{ij}=a_{ji}=0 \textrm{ or } -1$$
for $i\neq j \in \mathcal{V}$.

Then the Kac-Moody Lie algebra $\mathfrak{KM}(A)$ is the free Lie algebra over $\mathbb{F}$ with generators
$$e_i,f_i,h_i,\textrm{ where }i\in \mathcal{V}$$
subject to the relations
$$[h_i,h_j]=0,\ [h_i,e_j]=a_{ij}e_j,\ [h_i,f_j]=-a_{ij}f_j \  
\textrm{ for all } i,j\in \mathcal{V}$$
$$[e_i,f_j]=0,\ [e_i,f_i]=h_i,\ [e_i,[e_i,e_j]]=0=[f_i,[f_i,f_j]] \textrm{ for all } i\neq j \in \mathcal{V}.$$

The so-called compact subalgebra $\mathfrak{K}(A)$ of $\mathfrak{KM}(A)$ is the Lie subalgebra
generated by the elements

$$e_i+f_i, \ i\in \mathcal{V}.$$

If for each $x\in \mathcal{V}$ we denote by $x$ the element $e_x+f_x$, and consider the associated labeled graph $\Gamma=(\mathcal{V},\mathcal{E},\lambda)$ with vertex set $\mathcal{V}$ and two distinct vertices $x,y$ adjacent
if and only if $a_{xy}\neq 0$, and all labels $-1$, then we obtain the following.

\begin{lemma}
Let $x\neq y\in \mathcal{V}$. Then
$$\begin{array}{rll}
[x,y]&=0 & \textrm{ if }(x,y)\not\in \mathcal{E}\\
{[x,[x,y]]}&=-y & \textrm{ if }(x,y)\in \mathcal{E}.\\
\end{array}$$
\end{lemma}

\begin{proof}
If $x$ and $y$ are non-adjacent then clearly $[x,y]=0$. So, assume $x$ and $y$ are adjacent. Then

$$\begin{array}{ll}
[x,[x,y]]&=[e_x+f_x,[e_x+f_x,e_y+f_y]]]\\
&=[e_x+f_x,[e_x,e_y]+[f_x,f_y]]\\
&=[e_x,[e_x,e_y]]+ [e_x,[f_x,f_y]]+[f_x,[e_x,e_y]]+ [f_x,[f_x,f_y]]\\
&=[e_x,[f_x,f_y]]+[f_x,[e_x,e_y]]\\
&=-[f_y,[e_x,f_x]]-[e_y,[f_x,e_x]]\\
&=-[f_y,h_x]+[e_y,h_x]\\
&=-f_y-e_y\\
&=-y.
\end{array}$$
\end{proof}

\begin{theorem}(Berman \cite{Berman})\label{theorem_Berman}
Let $\mathbb{F}$ be a field of characteristic $0$ and $A=(a_{ij})$  a simply laced generalized Cartan matrix with associated graph $\Gamma=(\mathcal{V},\mathcal{E},\lambda)$.
Let $\mathfrak{KM}(A)$ be the Kac-Moody  Lie algebra over $\mathbb{F}$.

Then the compact Lie subalgebra 
$\mathfrak{K}(A)$ of $\mathfrak{KM}(A)$  is isomorphic to the free Lie algebra $\mathfrak{g}_\Gamma$ over $\mathbb{F}$ generated by $\mathcal{V}$
subject to the relations
$$\begin{array}{lll}
[x,y]&=0 & \textrm{ if }(x,y)\not\in \mathcal{E}\\
{[x,[x,y]]}&=-y & \textrm{ if }(x,y)\in \mathcal{E}\\
\end{array}$$
for $x\neq y\in \mathcal{V}$.
\end{theorem}

Combining the above \cref{theorem_Berman} with \cref{spin_exists}, we obtain the following.

\begin{corollary}\label{cor_Berman}
Let $\mathbb{F}$ be a field of characteristic $0$ and $A=(a_{ij})$  a simply laced generalized Cartan matrix with associated graph $\Gamma=(\mathcal{V},\mathcal{E},\lambda)$.
Let $\mathfrak{KM}(A)$ the Kac-Moody  Lie algebra over $\mathbb{F}$.

The compact Lie subalgebra 
$\mathfrak{K}(A)$ of $\mathfrak{KM}(A)$ admits a quotient isomorphic to $\mathfrak{K}(\Gamma)$, and in particular, admits a generalized spin representation.
\end{corollary}

\begin{example}
If $\Gamma$ is the graph $E_{10}$, with all vertices labeled by $-1$, and $\mathbb{F}$ is a field of type III, for example $\mathbb{R}$, then 
$\mathfrak{C}(\Gamma)$ is isomorphic to $\mathfrak{C}(V,Q)$, where $(V,Q)$ is a nondegenerate form of $+$-type.
But then $\mathfrak{g}(\Gamma)=\mathfrak{K}(\Gamma)$ is 
isomorphic to $\mathfrak{so}(32,\mathbb{F})$.
So, if $\mathbb{F}=\mathbb{R}$, we find that the compact Lie subalgebra 
$\mathfrak{K}(E_{10})$ of $\mathfrak{KM}(E_{10})$ admits a quotient isomorphic to $\mathfrak{so}(32,\mathbb{F})$.
Using \cref{EnCliffordTable}, we obtain similar results for graphs of type $E_n$ for all $n$. 
See also \cite{Buyl,Damour,spin}.
\end{example}

\begin{example}\label{compact}
As in \cite{spin} we can use the above result to determine the maximal compact Lie subalgebra $\mathfrak{K}$ of the semi-simple split real Lie algebras 
of type $A_n$, $D_n$, where $n\geq 1$ and $E_n$, where $6\leq n\leq 8$.

Indeed, using \ref{cor_Berman}, we find that the maximal compact Lie subalgebra  $\mathfrak{K}$ of a semi-simple split real Lie algebras $\mathfrak{g}$
of type $A_n$, $D_n$, where $n\geq 1$ and $E_n$, where $6\leq n\leq 8$, admits a quotient isomorphic to $\mathfrak{K}(\Gamma)$,
where $\Gamma$ is the corresponding graph of type $A_n$, $D_n$, or $E_n$ with all labels $-1$.

Using the results of the previous sections, we find  these quotients to be as in Table \ref{maxcompact}.
This provides a lowerbound for the dimension of   $\mathfrak{K}$ which coincides with the upperbound of the dimension of $\mathfrak{K}$ that one can obtain
from the  Iwasawa decomposition  of $\mathfrak{g}$.

These results can also be found in \cite{Tanya2}, where  $\mathfrak{K}$ is embedded in the Lie algebra of a (generalized) Clifford algebra.

\begin{table}
\begin{tabular}{|l|l|l|}\hline
Type of $\mathfrak{g}$ & Maximal compact subalgebra  $\mathfrak{K}$& dimension\\
\hline\hline
$A_n$  & $\mathfrak{so}(n+1,\mathbb{R})$& $\binom{n+1}{2}$\\
$D_n \ (n>3)$ & $\mathfrak{so}(n,\mathbb{R})\oplus \mathfrak{so}(n,\mathbb{R})$&$n(n-1)$\\
$E_6$       & $\mathfrak{sp}(4,\mathbb{H})$& $36$\\
$E_7$      & $\mathfrak{su}(8,\mathbb{C})$ & $63$\\
$E_8$       & $\mathfrak{so}(16,\mathbb{R})$ & $120$\\
\hline
\end{tabular}

\caption{Maximal compact subalgebras of the split real Lie algebras.}\label{maxcompact}

\end{table}
\end{example}

\section{Jordan Algebras}
\label{sect:jordan}

In this final section we consider Jordan algebras.
Let $\Gamma=(\mathcal{V},\mathcal{E},\lambda)$ be a graph, in which each vertex $v$ labeled with an invertible element $\lambda(v)$ of a ring $R$
containing distinct $-1$ and $1$ and consider $\mathfrak{C}:=\mathfrak{C}(\Gamma)$ which we identify with $\mathfrak{C}(V,g,\lambda)$, as in
the previous sections.

Then  $\mathfrak{C}^+$ is the Jordan algebra defined on   $\mathfrak{C}$ with the product $\circ$
defined as the linear extension of the product on the basis elements $v,w\in V$ given by
$$
v\circ w=\begin{cases}
0 &\mathrm{if}\ g(v,w)+g(w,v)=1\\
vw & \mathrm{if}\ g(v,w)+g(w,v)=0.\\
\end{cases}
$$

The elements $v,w\in \mathcal{V}$ then satisfy the following relations:

$$\begin{array}{ll}
v\circ(v\circ w)
&=\begin{cases}
0 &\mathrm{if} \ v\sim w\\
\lambda(v) w& \mathrm{if}\ v\not\sim w.\\
%= - \Lambda(v,v) w & \mathrm{if}\ g(v,w)+g(w,v)=0, g(v,v)=1.\\
%= \Lambda(v,v) w & \mathrm{if}\ g(v,w)=g(w,v)=0, g(v,v)=0\\
%= - \Lambda(v,v) w & \mathrm{if}\ g(v,w)=g(w,v)=1, g(v,v)=1\\
%=  \Lambda(v,v) w & \mathrm{if}\ g(v,w)=g(w,v)=1, g(v,v)=0\\
%= -\Lambda(v,v) w & \mathrm{if}\ g(v,w)=g(w,v)=0, g(v,v)=1\\
\end{cases}
\end{array}$$

By $\mathfrak{J}:=\mathfrak{J}(\Gamma)$ we denote the Jordan subalgebra of $\mathfrak{C}^+$ generated by the vertices of $\Gamma$.
If the complement of the graph $\Gamma$ is disconnected, we find that this subalgebra is a direct sum of Jordan subalgebras.
So, for the remainder of this section, assume that the complement of $\Gamma$ is connected.

There where the cotriangular space on $(V,f)$ was used to study the Lie algebra $\mathfrak{C}^-$ and its subalgebras,
we can now use the polar space $\mathcal{S}(V,f)$ on $(V,f)$ to study $\mathfrak{C}^+$ and its subalgebras.

Denote by $\mathcal{S}(\Gamma)$ the subspace of $\mathcal{S}(V,f)$ generated by the vertices of $\Gamma$.
Then we have:

\begin{proposition}\label{polarisbasis}
$\mathcal{S}(\Gamma)\cup \{\one\}$ is a basis for $\mathfrak{J}$.
\end{proposition}

\begin{proof}
If $v,w\in V$ are collinear points of $\mathcal{S}(\Gamma)$, then $v\circ w$ is a scalar multiple of  $v\pluscirc w$. So
$v\pluscirc w$ is also in $\mathcal{S}(\Gamma)$. So,  $\mathfrak{J}$ is linearly spanned by the elements from $\mathcal{S}(\Gamma)$ together with $\one$.
\end{proof}

\begin{proposition}\label{hyperplanepolar}
Suppose  $\Gamma$ and its complement are connected.

Then there is a unique quadratic form $Q$ compatible with $\Gamma$
such that  $\mathcal{S}(\Gamma)=\{v\in V\setminus\{0\}\mid  Q(v)=0\}$.
\end{proposition}

\begin{proof}
Let $\mathcal{R}$ be the radical of the polar space $\mathcal{S}(V,f)$ and let $\overline{\mathcal{S}}$ be the quotient of $\mathcal{S}(V,f)$ modulo this radical. Then  $\overline{\mathcal{S}}$ is isomorphic to $\mathcal{S}(\overline{V},\overline{f})$, where $\overline{V}=V/\Rad(f)$ and  $\overline{f}$ is the alternating form induced by $f$ on $\overline{V}$.

Suppose that $\overline{V}$ has dimension at least $6$. 
The subspace $\overline{\mathcal{S}_0}$ of $\overline{\mathcal{S}}$ generated by the points $p+\mathcal{R}$ with $p\in \mathcal{S}(\Gamma)$ is then equal to the   subspace 
$\mathcal{S}(\overline{V},\overline{Q})$, where $\overline{Q}$ is the quadratic form on $\overline{V}$ with associated bilinear form $\overline{f}$
which takes the value $0$ on the generators $v+\Rad(f)$, for all $v\in \mathcal{V}$.
This is well known, see \cite[Chapter 10]{BC}, or  the recent work  \cite{Pasini}.

If $\overline{V}$ has dimension $4$, then either $\overline{\mathcal{S}_0}$ is a grid with $9$ points, or $\overline{\mathcal{S}_0}$ is a set of at least $4$ pairwise noncollinear points, contained in the set of $5$ isotropic points with respect to $\overline{Q}$. 
In the first case we find $\overline{\mathcal{S}_0}$ to be the set of all isotropic points of $\overline{Q}$ (which is then of hyperbolic type).
In the second case $\overline{\mathcal{S}_0}$ consists of pairwise noncollinear  isotropic points of $\overline{Q}$, which is then of elliptic type.
This second case, however, can not occur, as  we find that either $\mathcal{V}$ contains points in the radical of ${\mathcal{S}}$ contradicting that
$\Gamma$ is connected, or none of the points of $\mathcal{V}$ is in the radical, contradicting the complement of $\Gamma$ to be connected.
We conclude that $\overline{\mathcal{S}_0}$ consists of all the isotropic points of $\overline{Q}$.

Next, we notice that the radical $\mathcal{R}$ is contained in $\mathcal{S}(\Gamma)$. Indeed, if this would not be the case, then the radical of 
$\mathcal{S}(\Gamma)$ is a proper subspace of $\mathcal{R}$ and the vectors in $\mathcal{V}$ would generate a proper subspace of $V$.

Combining this information with the observation that $\mathcal{S}(\Gamma)$ is contained in $\mathcal{S}(V,Q)$, where $Q$ is the unique quadratic form compatible with $\Gamma$ which takes the value $0$ on the vertices of $\Gamma$, we find that $\mathcal{S}(\Gamma)$ equals  $\mathcal{S}(V,Q)$.
\end{proof}

\begin{theorem}\label{jordancentralizerthm}
Suppose  $\Gamma$ and its complement are connected.

Then there is a  unique quadratic form $Q$ compatible with  $\Gamma$ and $0$ on the vertices of $\Gamma$  such that
$\mathfrak{J}$ is the centralizer of $\tau_Q$ in $\mathfrak{C}^+$.
\end{theorem}

\begin{proof}
This follows directly from \cref{hyperplanepolar}.
\end{proof}

We can now prove \cref{mainjordan}.
\begin{theorem}\label{mainjordansection}
Suppose  $\Gamma$ and its complement are connected and the ring $R$ contains $\frac{1}{2}$.

Let $\mathfrak{A}$ be an associative algebra containing  elements $\overline{v}$, with $v\in\mathcal{V}$, satisfying the following relations 
in $\mathfrak{A}^+$:
$$\begin{array}{lll}
   \overline{v}^2&= \lambda(v)\one_{\mathfrak{A}}, &  \\
   {\overline{v}\circ\overline{w}}&  =0 & \textrm{if } v\sim w,\\
   {\overline{v}\circ(\overline{v}\circ\overline{w})}&= \lambda(v) \overline{w} & \textrm{if } v \not\sim w.\\
\end{array}$$
Here $\circ$ denotes the Jordan product.

Then the subalgebra $\mathfrak{J}_0$ of  $\mathfrak{A}^+$ generated by the elements $\overline{v}$, with $v\in\mathcal{V}$, is isomorphic to 
a quotient of the subalgebra $\{x\in \mathfrak{C}^+\mid \tau_Q(x)=x\}$ of $\mathfrak{C}^+$ of symmetric elements with respect to the involution $\tau_{Q}$, where  $Q$ is the unique quadratic form compatible with $\Gamma$ and with $Q(v)=0$ for all $v\in \mathcal{V}$.
\end{theorem}

\begin{proof}
We  calculate in $\mathfrak{A}$.

Suppose $u,v\in \mathcal{V}$.
Then either $\overline{u}\overline{v}+\overline{v}\overline{u}=0$, or 
$$\begin{array}{ll}\overline{u}\circ(\overline{u}\circ \overline{v})&=\frac{1}{4}(\overline{u}(\overline{u}\overline{v}+\overline{v}\overline{u})+(\overline{u}\overline{v}+\overline{v}\overline{u})\overline{u})\\
&=\frac{1}{4}(\lambda(u)\overline{v}+2\overline{u}\overline{v}\overline{u}+\lambda(u)\overline{v})\\
&=\frac{1}{2}(\lambda(u)\overline{v}+\overline{u}\overline{v}\overline{u})\\
&=\lambda(u) \overline{v}.\\
\end{array}$$
So, $\lambda(u)\overline{v}-\overline{u}\overline{v}\overline{u}=0$. Multiplying with $\overline{u}$ yields
$\lambda(u)(\overline{u}\overline{v}-\overline{v}\overline{u})=0$ 
and thus $\overline{u}\overline{v}-\overline{v}\overline{u}=0$.

So (the images of) the elements in $\mathcal{V}$ satisfy the relations of the generators of $\mathfrak{C}:=\mathfrak{C}(\Gamma)$.
In particular, they generate a subalgebra of  $\mathfrak{A}$, which is isomorphic to a quotient of $\mathfrak{C}$.
This implies that $\mathfrak{J}_0$ is isomorphic to a quotient of the subalgebra $\mathfrak{J}(\Gamma)$ of $\mathfrak{C}^+$, which, by \cref{jordancentralizerthm}, is the algebra
 $\{x\in \mathfrak{C}^+\mid \tau_Q(x)=x\}$ of symmetric elements with respect to the involution $\tau_{Q}$, where  $Q$ is the unique quadratic form compatible with $\Gamma$ and with $Q(v)=0$ for all $v\in \mathcal{V}$.
\end{proof}

\begin{remark}
By  famous results of Zelmanov (see \cite[Section 7.3]{mccrimmon}),  simple Jordan algebras over a field $\mathbb{F}$ of characteristic $\neq 2$ are special, except maybe when they have dimension $27$ over their center.
So, a simple Jordan algebra over such field $\mathbb{F}$ generated by elements $\overline{v}$, where $v\in \mathcal{V}$ the vertex set of a graph $\Gamma$
which is connected and coconnected, satisfying the relations
$$\begin{array}{lll}
   \overline{v}^2&= \lambda(v)\one_{\mathfrak{A}}, &  \\
   {\overline{v}\circ\overline{w}}&  =0 & \textrm{if } v\sim w,\\
   {\overline{v}\circ(\overline{v}\circ\overline{w})}&= \lambda(v) \overline{w} & \textrm{if } v \not\sim w.\\
\end{array}$$
is either of dimension $27$ or does satisfy the conclusion of the above theorem.
\end{remark}

 We end this paper with an example.

 \begin{example}
 Let $\Gamma$ be the $E_6$ graph with all vertices labeled  $-1$.
 Then consider the special real Jordan algebra $\mathfrak{J}(\Gamma)$
 generated by the $6$ vertices of $\Gamma$ with relations
 $$\begin{array}{lll}
 v^2=v\circ v&=-\one\\
 v\circ w&=0&\mathrm{for} \  v\sim w\\
 v\circ(v\circ w)&= -w& \mathrm{for}\ v\not\sim w\\
 \end{array}$$
 where $v,w\in \mathcal{V}$.
 
 The quadratic form $Q$ compatible with $\Gamma$ having $Q(v)$ equal to $1$ for all $v\in \mathcal{V}$ is of $-$-type.
 So, the algebra $\mathfrak{C}=\mathfrak{C}(\Gamma)$ is the algebra of $4\times 4$ matrices with coefficients in the 
 quaternions $\mathbb{H}$, see also \cref{EnCliffordTable}.
 
 The involution $\tau_{Q}$ induces the map 
 $$M\mapsto \overline{M}^\top,$$ where $\overline{\phantom{a}}$ denotes quaternion conjugation.
 
 The quadratic form $Q'$ compatible with $\Gamma$ having $Q'(v)=0$ for all $v\in \mathcal{V}$, however, is of $+$-type.
 (Indeed, in the diagram one finds already $3$  nodes generating an isotropic $3$-space of $V$.)
 As there are $35$ isotropic vectors with respect to $Q'$, we find 
 the Jordan subalgebra $\mathfrak{J}_0$ of $\mathfrak{C}^+$ generated by the vertices of $\Gamma$  
 to be $36$-dimensional ($\one$ is also in $\mathfrak{J}_0$). 
 The anti-automorphism $\tau_{Q'}$ on $\mathfrak{C}$ corresponds to the map
 $$M\mapsto -J\overline{M}^\top J,$$
 where $J=\begin{pmatrix}0&1\\-1&0\end{pmatrix}$.
 So  $\mathfrak{J}_0$ is isomorphic to the special Jordan algebra of $4\times 4$ matrices $\begin{pmatrix}A&B\\C&D\\\end{pmatrix}$ with coefficients in the quaternions 
 $\mathbb{H}$ satisfying 
 $$\begin{pmatrix}A&B\\C&D\\\end{pmatrix}=\begin{pmatrix}\overline{D}^\top&-\overline{B}^\top\\-\overline{C}^\top&\overline{A}^\top\\\end{pmatrix}.$$

 The above results show that $\mathfrak{J}$ is isomorphic to a quotient of $\mathfrak{J}_0$.
 As the latter is simple, see \cref{simple}, we find $\mathfrak{J}$ to be isomorphic to $\mathfrak{J}_0$. 
 \end{example}

\newpage

\bibliographystyle{plain}

\bibliography{clifford.bib}

@article {Hallembed,
    AUTHOR = {Hall, Jonathan I.},
     TITLE = {Linear representations of cotriangular spaces},
   JOURNAL = {Linear Algebra Appl.},
  FJOURNAL = {Linear Algebra and its Applications},
    VOLUME = {49},
      YEAR = {1983},
     PAGES = {257--273},
      ISSN = {0024-3795},
   MRCLASS = {51E20 (05B30 51F25 94B05)},
  MRNUMBER = {688390},
MRREVIEWER = {Joseph A. Thas},
       DOI = {10.1016/0024-3795(83)90108-8},
       URL = {https://doi.org/10.1016/0024-3795(83)90108-8},
}

@article {leopardi1,
    AUTHOR = {Leopardi, Paul},
     TITLE = {Twin bent functions, strongly regular {C}ayley graphs, and
              {H}urwitz-{R}adon theory},
   JOURNAL = {J. Algebra Comb. Discrete Struct. Appl.},
  FJOURNAL = {Journal of Algebra Combinatorics Discrete Structures and
              Applications},
    VOLUME = {4},
      YEAR = {2017},
    NUMBER = {3},
     PAGES = {271--280},
   MRCLASS = {05E30 (06E30 11T71 15B34)},
  MRNUMBER = {3673231},
MRREVIEWER = {Ilia N. Ponomarenko},
       DOI = {10.13069/jacodesmath.327377},
       URL = {https://doi.org/10.13069/jacodesmath.327377},
}

@article {Seberry1,
    AUTHOR = {Seberry, Jennifer},
     TITLE = {Hadamard matrices, orthogonal designs and
              {C}lifford--{G}astineau-{H}ills algebras},
   JOURNAL = {Australas. J. Combin.},
  FJOURNAL = {The Australasian Journal of Combinatorics},
    VOLUME = {71},
      YEAR = {2018},
     PAGES = {452--467},
      ISSN = {1034-4942},
   MRCLASS = {05B15 (05B20)},
  MRNUMBER = {3801276},
MRREVIEWER = {Jianying Fang},
}

@article {Leopardi2,
    AUTHOR = {Leopardi, Paul C.},
     TITLE = {Constructions for {H}adamard matrices using {C}lifford
              algebras, and their relation to amicability/anti-amicability
              graphs},
   JOURNAL = {Australas. J. Combin.},
  FJOURNAL = {The Australasian Journal of Combinatorics},
    VOLUME = {58},
      YEAR = {2014},
     PAGES = {214--248},
      ISSN = {1034-4942},
   MRCLASS = {15B05 (05B20 15A66)},
  MRNUMBER = {3211780},
MRREVIEWER = {Zhanna Charniak},
}

@incollection {Kaplansky,
    AUTHOR = {Kaplansky, Irving},
     TITLE = {Some simple {L}ie algebras of characteristic {$2$}},
 BOOKTITLE = {Lie algebras and related topics ({N}ew {B}runswick, {N}.{J}.,
              1981)},
    SERIES = {Lecture Notes in Math.},
    VOLUME = {933},
     PAGES = {127--129},
 PUBLISHER = {Springer, Berlin-New York},
      YEAR = {1982},
   MRCLASS = {17B20},
  MRNUMBER = {675111},
}

@book {BC,
   author={Buekenhout, F. and Cohen, A. M.},
   title={Diagram geometry, Related to classical groups and buildings},
   series={Ergebnisse der Mathematik und ihrer Grenzgebiete. 3. Folge. A
   Series of Modern Surveys in Mathematics },
   volume={57},
   publisher={Springer, Heidelberg},
   date={2013},
   pages={xiv+592},
   isbn={978-3-642-34452-7},
   isbn={978-3-642-34453-4},
   year={2013},
}

@article {Erik,
    AUTHOR = {Brouwer, Andries E. and Cohen, Arjeh M. and Cuypers, Hans and
              Hall, Jonathan I. and Postma, Erik},
     TITLE = {Lie algebras, 2-groups and cotriangular spaces},
   JOURNAL = {Adv. Geom.},
  FJOURNAL = {Advances in Geometry},
    VOLUME = {12},
      YEAR = {2012},
    NUMBER = {1},
     PAGES = {1--17},
      ISSN = {1615-715X},
   MRCLASS = {51A50 (17B99)},
  MRNUMBER = {2911156},
MRREVIEWER = {Harald L\"{o}we},
       DOI = {10.1515/advgeom.2011.047},
       URL = {https://doi.org/10.1515/advgeom.2011.047},
}

@article {cotriangleLie,
    AUTHOR = {Cuypers, Hans},
     TITLE = {Lie algebras and cotriangular spaces},
   JOURNAL = {Bull. Belg. Math. Soc. Simon Stevin},
  FJOURNAL = {Bulletin of the Belgian Mathematical Society. Simon Stevin},
    VOLUME = {12},
      YEAR = {2005},
    NUMBER = {2},
     PAGES = {209--221},
      ISSN = {1370-1444},
   MRCLASS = {17B40 (51E26)},
  MRNUMBER = {2179964},
MRREVIEWER = {Jonathan I. Hall},
       URL = {http://projecteuclid.org/euclid.bbms/1117805084},
}

@article {quasi,
    AUTHOR = {Gastineau-Hills, H. M.},
     TITLE = {Quasi-{C}lifford algebras and systems of orthogonal designs},
   JOURNAL = {J. Austral. Math. Soc. Ser. A},
  FJOURNAL = {Australian Mathematical Society. Journal. Series A. Pure
              Mathematics and Statistics},
    VOLUME = {32},
      YEAR = {1982},
    NUMBER = {1},
     PAGES = {1--23},
      ISSN = {0263-6115},
   MRCLASS = {15A66 (05B15)},
  MRNUMBER = {643426},
MRREVIEWER = {Maurice Flamant},
}

@article {Elduque_clifford,
    AUTHOR = {Elduque, Alberto and Rodrigo-Escudero, Adri\'{a}n},
     TITLE = {Clifford algebras as twisted group algebras and the {A}rf
              invariant},
   JOURNAL = {Adv. Appl. Clifford Algebras},
  FJOURNAL = {Advances in Applied Clifford Algebras},
    VOLUME = {28},
      YEAR = {2018},
    NUMBER = {2},
     PAGES = {Art. 41, 15},
      ISSN = {0188-7009},
   MRCLASS = {15A66 (11E04 15A63 16S35 16W50)},
  MRNUMBER = {3790094},
MRREVIEWER = {Alberto Molgado},
       DOI = {10.1007/s00006-018-0862-y},
       URL = {https://doi.org/10.1007/s00006-018-0862-y},
}

@article {Twisted_clifford,
    AUTHOR = {Albuquerque, Helena and Majid, Shahn},
     TITLE = {Clifford algebras obtained by twisting of group algebras},
   JOURNAL = {J. Pure Appl. Algebra},
  FJOURNAL = {Journal of Pure and Applied Algebra},
    VOLUME = {171},
      YEAR = {2002},
    NUMBER = {2-3},
     PAGES = {133--148},
      ISSN = {0022-4049},
   MRCLASS = {15A66 (16W30)},
  MRNUMBER = {1904473},
       DOI = {10.1016/S0022-4049(01)00124-4},
       URL = {https://doi.org/10.1016/S0022-4049(01)00124-4},
}

@article {Damour,
    AUTHOR = {Damour, Thibault and Kleinschmidt, Axel and Nicolai, Hermann},
     TITLE = {Hidden symmetries and the fermionic sector of
              eleven-dimensional supergravity},
   JOURNAL = {Phys. Lett. B},
  FJOURNAL = {Physics Letters. B. Particle Physics, Nuclear Physics and
              Cosmology},
    VOLUME = {634},
      YEAR = {2006},
    NUMBER = {2-3},
     PAGES = {319--324},
      ISSN = {0370-2693},
   MRCLASS = {83E50 (81T30 83E30)},
  MRNUMBER = {2202945},
       DOI = {10.1016/j.physletb.2006.01.015},
       URL = {https://doi.org/10.1016/j.physletb.2006.01.015},
}

@article {Buyl,
    AUTHOR = {de Buyl, Sophie and Henneaux, Marc and Paulot, Louis},
     TITLE = {Extended {$E_8$} invariance of 11-dimensional supergravity},
   JOURNAL = {J. High Energy Phys.},
  FJOURNAL = {Journal of High Energy Physics. A SISSA Journal},
      YEAR = {2006},
    NUMBER = {2},
     PAGES = {056, 11},
      ISSN = {1126-6708},
   MRCLASS = {83E50 (81T30 83E30)},
  MRNUMBER = {2219449},
MRREVIEWER = {Luis Joaqu\'{\i}n Boya},
       DOI = {10.1088/1126-6708/2006/02/056},
       URL = {https://doi.org/10.1088/1126-6708/2006/02/056},
    volume = {2006},
}

@book {Seidel,
    AUTHOR = {Seidel, J. J.},
     TITLE = {On two-graphs and {S}hult's characterization of symplectic and
              orthogonal geometries over {GF}(2)},
      NOTE = {T.H.-Report, No. 73-WSK-02},
 PUBLISHER = {Department of Mathematics, Technological University Eindhoven,
              Eindhoven},
      YEAR = {1973},
     PAGES = {i+25},
   MRCLASS = {05B25},
  MRNUMBER = {0422055},
MRREVIEWER = {Francis Buckenhout},
}

@article {Shult,
    AUTHOR = {Shult, Ernest E.},
     TITLE = {Characterizations of certain classes of graphs},
   JOURNAL = {J. Combinatorial Theory Ser. B},
  FJOURNAL = {Journal of Combinatorial Theory. Series B},
    VOLUME = {13},
      YEAR = {1972},
     PAGES = {142--167},
      ISSN = {0095-8956},
   MRCLASS = {05C99 (20F99)},
  MRNUMBER = {311518},
MRREVIEWER = {J. J. Seidel},
       DOI = {10.1016/0095-8956(72)90050-0},
       URL = {https://doi.org/10.1016/0095-8956(72)90050-0},
}

@article {Hall,
    AUTHOR = {Hall, J. I.},
     TITLE = {Graphs, geometry, {$3$}-transpositions, and symplectic {${\bf
              F}_2$}-transvection groups},
   JOURNAL = {Proc. London Math. Soc. (3)},
  FJOURNAL = {Proceedings of the London Mathematical Society. Third Series},
    VOLUME = {58},
      YEAR = {1989},
    NUMBER = {1},
     PAGES = {89--111},
      ISSN = {0024-6115},
   MRCLASS = {20F32 (05C75)},
  MRNUMBER = {969549},
       DOI = {10.1112/plms/s3-58.1.89},
       URL = {https://doi.org/10.1112/plms/s3-58.1.89},
}

@article {Berman,
    AUTHOR = {Berman, S.},
     TITLE = {On generators and relations for certain involutory subalgebras
              of {K}ac-{M}oody {L}ie algebras},
   JOURNAL = {Comm. Algebra},
  FJOURNAL = {Communications in Algebra},
    VOLUME = {17},
      YEAR = {1989},
    NUMBER = {12},
     PAGES = {3165--3185},
      ISSN = {0092-7872},
   MRCLASS = {17B67},
  MRNUMBER = {1030614},
MRREVIEWER = {Sen Qiu},
       DOI = {10.1080/00927878908823899},
       URL = {https://doi.org/10.1080/00927878908823899},
}

@book {Seberry,
    AUTHOR = {Seberry, Jennifer},
     TITLE = {Orthogonal designs},
      NOTE = {Hadamard matrices, quadratic forms and algebras,
              Revised and updated edition of the 1979 original [ MR0534614]},
 PUBLISHER = {Springer, Cham},
      YEAR = {2017},
     PAGES = {xxiii+453},
      ISBN = {978-3-319-59031-8; 978-3-319-59032-5},
   MRCLASS = {05-02 (05B05 05B20 05B30)},
  MRNUMBER = {3729402},
MRREVIEWER = {Anthony B. Evans},
       DOI = {10.1007/978-3-319-59032-5},
       URL = {https://doi.org/10.1007/978-3-319-59032-5},
}

@article {Hall_extra,
    AUTHOR = {Hall, J. I.},
     TITLE = {The number of trace-valued forms and extraspecial groups},
   JOURNAL = {J. London Math. Soc. (2)},
  FJOURNAL = {Journal of the London Mathematical Society. Second Series},
    VOLUME = {37},
      YEAR = {1988},
    NUMBER = {1},
     PAGES = {1--13},
      ISSN = {0024-6107},
   MRCLASS = {03C60 (20F99)},
  MRNUMBER = {921742},
       DOI = {10.1112/jlms/s2-37.121.1},
       URL = {https://doi.org/10.1112/jlms/s2-37.121.1},
}

@article {Shirokov18,
    AUTHOR = {Shirokov, D. S.},
     TITLE = {Classification of {L}ie algebras of specific type in
              complexified {C}lifford algebras},
   JOURNAL = {Linear Multilinear Algebra},
  FJOURNAL = {Linear and Multilinear Algebra},
    VOLUME = {66},
      YEAR = {2018},
    NUMBER = {9},
     PAGES = {1870--1887},
      ISSN = {0308-1087},
   MRCLASS = {17B60 (15A66 22E60)},
  MRNUMBER = {3825703},
MRREVIEWER = {P. D. Beites},
       DOI = {10.1080/03081087.2017.1376612},
       URL = {https://doi.org/10.1080/03081087.2017.1376612},
}

@article {Shirokov16,
    AUTHOR = {Shirokov, Dmitry},
     TITLE = {On some {L}ie groups containing spin group in {C}lifford
              algebra},
   JOURNAL = {J. Geom. Symmetry Phys.},
  FJOURNAL = {Journal of Geometry and Symmetry in Physics},
    VOLUME = {42},
      YEAR = {2016},
     PAGES = {73--94},
      ISSN = {1312-5192},
   MRCLASS = {22E46 (15A66)},
  MRNUMBER = {3586444},
       DOI = {10.7546/jgsp-42-2016-73-94},
       URL = {https://doi.org/10.7546/jgsp-42-2016-73-94},
}

@article {Shirokov15,
    AUTHOR = {Shirokov, D. S.},
     TITLE = {Symplectic, orthogonal and linear {L}ie groups in {C}lifford
              algebra},
   JOURNAL = {Adv. Appl. Clifford Algebras},
  FJOURNAL = {Advances in Applied Clifford Algebras},
    VOLUME = {25},
      YEAR = {2015},
    NUMBER = {3},
     PAGES = {707--718},
      ISSN = {0188-7009},
   MRCLASS = {22E46 (15A66)},
  MRNUMBER = {3384856},
MRREVIEWER = {Stephen Slebarski},
       DOI = {10.1007/s00006-014-0520-y},
       URL = {https://doi.org/10.1007/s00006-014-0520-y},
}

@article {seven,
    AUTHOR = {Seven, Ahmet I.},
     TITLE = {Orbits of groups generated by transvections over {${\Bbb
              F}_2$}},
   JOURNAL = {J. Algebraic Combin.},
  FJOURNAL = {Journal of Algebraic Combinatorics. An International Journal},
    VOLUME = {21},
      YEAR = {2005},
    NUMBER = {4},
     PAGES = {449--474},
      ISSN = {0925-9899},
   MRCLASS = {14L30 (05E15)},
  MRNUMBER = {2153935},
MRREVIEWER = {Shreedhar P. Inamdar},
       DOI = {10.1007/s10801-005-3021-1},
       URL = {https://doi.org/10.1007/s10801-005-3021-1},
}

@incollection {shaw1,
    AUTHOR = {Shaw, R.},
     TITLE = {Finite geometry, {D}irac groups and the table of real
              {C}lifford algebras},
 BOOKTITLE = {Clifford algebras and spinor structures},
    SERIES = {Math. Appl.},
    VOLUME = {321},
     PAGES = {59--99},
 PUBLISHER = {Kluwer Acad. Publ., Dordrecht},
      YEAR = {1995},
   MRCLASS = {15A66 (15A63)},
  MRNUMBER = {1483083},
MRREVIEWER = {Bernd A. Schmeikal},
}

@article {shaw2,
    AUTHOR = {Shaw, R.},
     TITLE = {Finite geometries and {C}lifford algebras},
   JOURNAL = {J. Math. Phys.},
  FJOURNAL = {Journal of Mathematical Physics},
    VOLUME = {30},
      YEAR = {1989},
    NUMBER = {9},
     PAGES = {1971--1984},
      ISSN = {0022-2488},
   MRCLASS = {15A66 (05B25)},
  MRNUMBER = {1009905},
MRREVIEWER = {Horia Scutaru},
       DOI = {10.1063/1.528234},
       URL = {https://doi.org/10.1063/1.528234},
}

@article {shaw3,
    AUTHOR = {Shaw, R. and Jarvis, T. M.},
     TITLE = {Finite geometries and {C}lifford algebras. {II}},
   JOURNAL = {J. Math. Phys.},
  FJOURNAL = {Journal of Mathematical Physics},
    VOLUME = {31},
      YEAR = {1990},
    NUMBER = {6},
     PAGES = {1315--1324},
      ISSN = {0022-2488},
   MRCLASS = {15A66},
  MRNUMBER = {1054318},
       DOI = {10.1063/1.528769},
       URL = {https://doi.org/10.1063/1.528769},
}

@incollection {shaw4,
    AUTHOR = {Shaw, R.},
     TITLE = {Finite geometries and {C}lifford algebras. {III}},
 BOOKTITLE = {Clifford algebras and their applications in mathematical
              physics ({M}ontpellier, 1989)},
    SERIES = {Fund. Theories Phys.},
    VOLUME = {47},
     PAGES = {121--132},
 PUBLISHER = {Kluwer Acad. Publ., Dordrecht},
      YEAR = {1992},
   MRCLASS = {15A66 (51E30)},
  MRNUMBER = {1199589},
MRREVIEWER = {Pierre Angl\`es},
}

@article {Beineke,
    AUTHOR = {Beineke, Lowell W.},
     TITLE = {Characterizations of derived graphs},
   JOURNAL = {J. Combinatorial Theory},
  FJOURNAL = {Journal of Combinatorial Theory},
    VOLUME = {9},
      YEAR = {1970},
     PAGES = {129--135},
      ISSN = {0021-9800},
   MRCLASS = {05.40},
  MRNUMBER = {262097},
}

@article {Tanya1,
    AUTHOR = {Khovanova, Tanya},
     TITLE = {Clifford algebras and graphs},
   JOURNAL = {Geombinatorics},
  FJOURNAL = {Geombinatorics},
    VOLUME = {20},
      YEAR = {2010},
    NUMBER = {2},
     PAGES = {56--76},
      ISSN = {1065-7371},
   MRCLASS = {15A66 (05C25)},
  MRNUMBER = {2732517},
}

@article {Tanya2,
    AUTHOR = {Khovanova, T. G.},
     TITLE = {Models of representations and generalized {C}lifford algebras},
   JOURNAL = {Funktsional. Anal. i Prilozhen.},
  FJOURNAL = {Akademiya Nauk SSSR. Funktsional\cprime ny\u{\i} Analiz i ego
              Prilozheniya},
    VOLUME = {16},
      YEAR = {1982},
    NUMBER = {4},
     PAGES = {90--91},
      ISSN = {0374-1990},
   MRCLASS = {22E46 (15A66)},
  MRNUMBER = {684144},
MRREVIEWER = {James E. Humphreys},
}

@article {gintz,
    AUTHOR = {Gintz, M.},
     TITLE = {Classifying Graph Lie Algebras},
   JOURNAL = {Preprint\\ https://math.mit.edu/research/highschool/primes/materials/2017/Gintz.pdf},
   year= {2018},
}

@article {alglinegraph,
    AUTHOR = {Roussopoulos, Nicholas D.},
     TITLE = {A max {$\{m,n\}$} algorithm for determining the graph {$H$}
              from its line graph {$G$}},
   JOURNAL = {Information Processing Lett.},
  FJOURNAL = {Information Processing Letters},
    VOLUME = {2},
      YEAR = {1973},
     PAGES = {108--112},
      ISSN = {0020-0190},
   MRCLASS = {94A20 (05C99)},
  MRNUMBER = {424435},
MRREVIEWER = {S. L. Hakimi},
       DOI = {10.1016/0020-0190(73)90029-x},
       URL = {https://doi.org/10.1016/0020-0190(73)90029-x},
}

@incollection {partition,
    AUTHOR = {Habib, Michel and Paul, Christophe and Viennot, Laurent},
     TITLE = {A synthesis on partition refinement: a useful routine for
              strings, graphs, {B}oolean matrices and automata},
 BOOKTITLE = {S{TACS} 98 ({P}aris, 1998)},
    SERIES = {Lecture Notes in Comput. Sci.},
    VOLUME = {1373},
     PAGES = {25--38},
 PUBLISHER = {Springer, Berlin},
      YEAR = {1998},
   MRCLASS = {68Q20 (68P05 68Q25)},
  MRNUMBER = {1650757},
       DOI = {10.1007/BFb0028546},
       URL = {https://doi.org/10.1007/BFb0028546},
}

@article {Ovsienko,
    AUTHOR = {Ovsienko, Valentin},
     TITLE = {Real {C}lifford algebras and quadratic forms over {$\Bbb
              F_2$}: two old problems become one},
   JOURNAL = {Math. Intelligencer},
  FJOURNAL = {The Mathematical Intelligencer},
    VOLUME = {38},
      YEAR = {2016},
    NUMBER = {3},
     PAGES = {1--5},
      ISSN = {0343-6993},
   MRCLASS = {11E88 (11E04 15A66)},
  MRNUMBER = {3546474},
MRREVIEWER = {Andrzej S\l adek},
       DOI = {10.1007/s00283-016-9637-2},
       URL = {https://doi.org/10.1007/s00283-016-9637-2},
}

@article {Pasini,
    AUTHOR = {Cardinali, I. and Giuzzi, L. and Pasini, A.},
     TITLE = {Nearly all subspaces of a classical polar space arise from its
              universal embedding},
   JOURNAL = {Linear Algebra Appl.},
  FJOURNAL = {Linear Algebra and its Applications},
    VOLUME = {627},
      YEAR = {2021},
     PAGES = {287--307},
      ISSN = {0024-3795},
   MRCLASS = {51A50 (51B25 51E24)},
  MRNUMBER = {4283256},
       DOI = {10.1016/j.laa.2021.06.013},
       URL = {https://doi.org/10.1016/j.laa.2021.06.013},
}

@article {spin,
    AUTHOR = {Hainke, Guntram and K\"{o}hl, Ralf and Levy, Paul},
     TITLE = {Generalized spin representations},
      NOTE = {With an appendix by Max Horn and Ralf K\"{o}hl},
   JOURNAL = {M\"{u}nster J. Math.},
  FJOURNAL = {M\"{u}nster Journal of Mathematics},
    VOLUME = {8},
      YEAR = {2015},
    NUMBER = {1},
     PAGES = {181--210},
      ISSN = {1867-5778},
   MRCLASS = {17B67 (17B10 17B25)},
  MRNUMBER = {3549525},
MRREVIEWER = {Stefano Capparelli},
       DOI = {10.17879/65219674985},
       URL = {https://doi.org/10.17879/65219674985},
}

@article {bouarroudj,
    AUTHOR = {Bouarroudj, Sofiane and Lebedev, Alexei and Leites, Dimitry
              and Shchepochkina, Irina},
     TITLE = {Lie algebra deformations in characteristic 2},
   JOURNAL = {Math. Res. Lett.},
  FJOURNAL = {Mathematical Research Letters},
    VOLUME = {22},
      YEAR = {2015},
    NUMBER = {2},
     PAGES = {353--402},
      ISSN = {1073-2780},
   MRCLASS = {17B50 (17B56)},
  MRNUMBER = {3342238},
MRREVIEWER = {Christopher M. Drupieski},
       DOI = {10.4310/MRL.2015.v22.n2.a3},
       URL = {https://doi.org/10.4310/MRL.2015.v22.n2.a3},
}

@book {mccrimmon,
    AUTHOR = {McCrimmon, Kevin},
     TITLE = {A taste of {J}ordan algebras},
    SERIES = {Universitext},
 PUBLISHER = {Springer-Verlag, New York},
      YEAR = {2004},
     PAGES = {xxvi+562},
      ISBN = {0-387-95447-3},
   MRCLASS = {17-02 (17Cxx)},
  MRNUMBER = {2014924},
MRREVIEWER = {Plamen Koshlukov},
}

\vspace{2cm}

\parindent=0pt
Hans Cuypers\\
Department of Mathematics and Computer Science \\
Eindhoven University of Technology\\
P.O. Box 513 5600 MB, Eindhoven\\
The Netherlands\\
email: f.g.m.t.cuypers@tue.nl

\end{document}